\documentclass[11pt]{amsart}
\usepackage{amsmath,amsthm,epsfig,amssymb}
\title{Half domination arrangements in regular and semi-regular tessellation type graphs}
\author{Eugen J. Ionascu}
\curraddr{Department of Mathematics\\ Columbus State University\\4225 University Avenue\\
Columbus, GA 31907\\
Honorific Member of the Romanian Institute of Mathematics ``Simion
Stoilow" } \email{ionascu@columbusstate.edu;}
\subjclass[2000]{ 05C35,\ 05C69}
\date{January $19^{th}$, 2012}
\begin{document}
\def\sms{\small\scshape}
\baselineskip18pt
\newtheorem{theorem}{\hspace{\parindent}
T{\scriptsize HEOREM}}[section]
\newtheorem{proposition}[theorem]
{\hspace{\parindent }P{\scriptsize ROPOSITION}}
\newtheorem{corollary}[theorem]
{\hspace{\parindent }C{\scriptsize OROLLARY}}
\newtheorem{lemma}[theorem]
{\hspace{\parindent }L{\scriptsize EMMA}}
\newtheorem{definition}[theorem]
{\hspace{\parindent }D{\scriptsize EFINITION}}
\newtheorem{problem}[theorem]
{\hspace{\parindent }P{\scriptsize ROBLEM}}
\newtheorem{conjecture}[theorem]
{\hspace{\parindent }C{\scriptsize ONJECTURE}}
\newtheorem{example}[theorem]
{\hspace{\parindent }E{\scriptsize XAMPLE}}
\newtheorem{remark}[theorem]
{\hspace{\parindent }R{\scriptsize EMARK}}
\renewcommand{\thetheorem}{\arabic{section}.\arabic{theorem}}
\renewcommand{\theenumi}{(\roman{enumi})}
\renewcommand{\labelenumi}{\theenumi}
\newcommand{\Q}{{\mathbb Q}}
\newcommand{\Z}{{\mathbb Z}}
\newcommand{\N}{{\mathbb N}}
\newcommand{\C}{{\mathbb C}}
\newcommand{\R}{{\mathbb R}}
\newcommand{\F}{{\mathbb F}}
\newcommand{\K}{{\mathbb K}}
\newcommand{\D}{{\mathbb D}}
\def\phi{\varphi}
\def\ra{\rightarrow}
\def\sd{\bigtriangledown}
\def\ac{\mathaccent94}
\def\wi{\sim}
\def\wt{\widetilde}
\def\bb#1{{\Bbb#1}}
\def\bs{\backslash}
\def\cal{\mathcal}
\def\ca#1{{\cal#1}}
\def\Bbb#1{\bf#1}
\def\blacksquare{{\ \vrule height7pt width7pt depth0pt}}
\def\bsq{\blacksquare}
\def\proof{\hspace{\parindent}{P{\scriptsize ROOF}}}
\def\pofthe{P{\scriptsize ROOF OF}
T{\scriptsize HEOREM}\  }
\def\pofle{\hspace{\parindent}P{\scriptsize ROOF OF}
L{\scriptsize EMMA}\  }
\def\pofcor{\hspace{\parindent}P{\scriptsize ROOF OF}
C{\scriptsize ROLLARY}\  }
\def\pofpro{\hspace{\parindent}P{\scriptsize ROOF OF}
P{\scriptsize ROPOSITION}\  }
\def\n{\noindent}
\def\wh{\widehat}
\def\eproof{$\hfill\bsq$\par}
\def\ds{\displaystyle}
\def\du{\overset{\text {\bf .}}{\cup}}
\def\Du{\overset{\text {\bf .}}{\bigcup}}
\def\b{$\blacklozenge$}
\def\eqtr{{\cal E}{\cal T}(\Z) }
\def\eproofi{\bsq}

\begin{abstract} We study the problem of half-domination sets of vertices in transitive infinite graphs
generated by regular or semi-regular tessellations of the plane. In some cases, the results obtained are
sharp and in the rest, we show upper bounds for the average densities of vertices in half-domination sets.
\end{abstract} \maketitle

\section{\label{intro}Introduction} 
\begin{figure}
$$\begin{array}{c}
\underset{(4^4)}{\epsfig{file=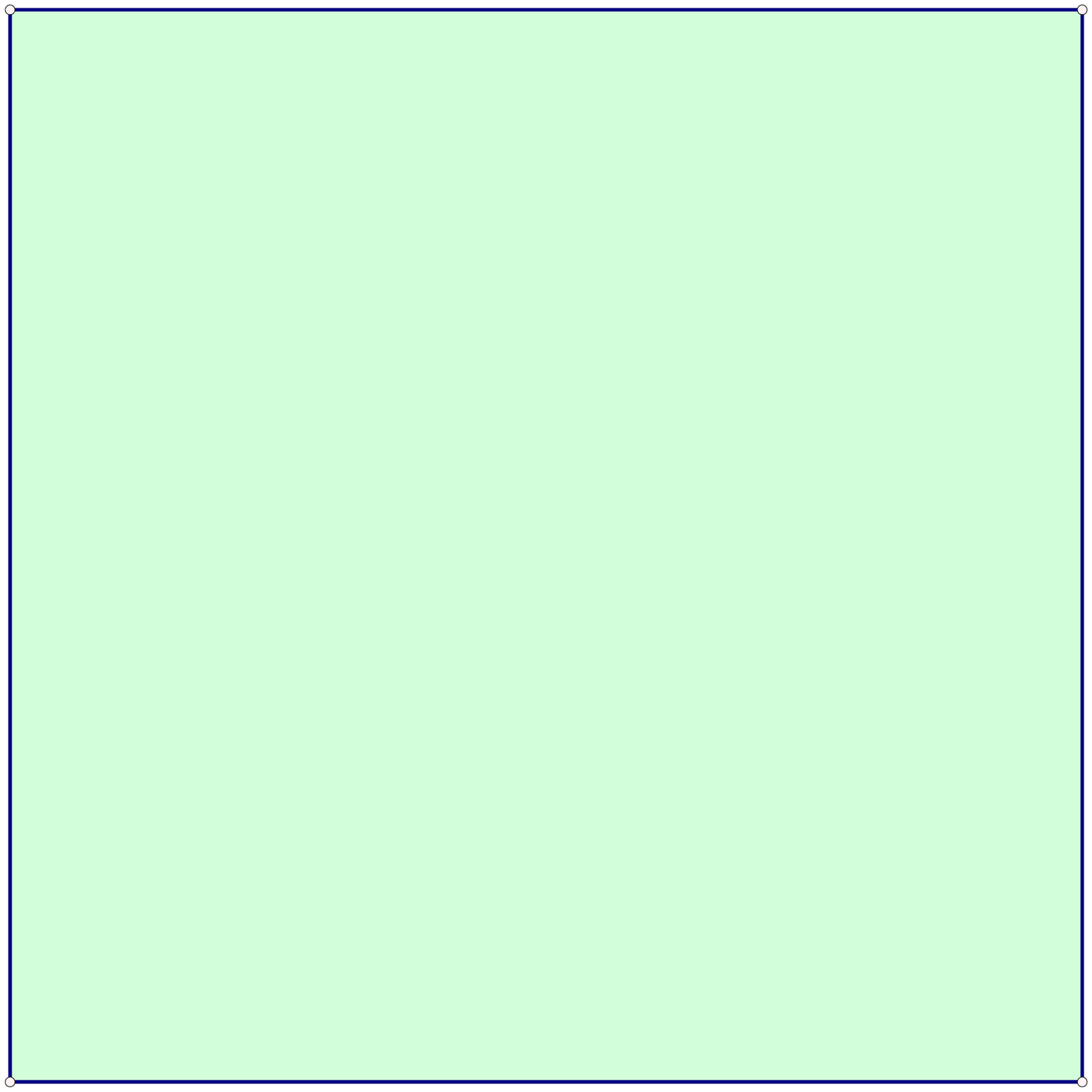,height=.5in,width=.5in}},\
\underset{(6^3)}{ \epsfig{file=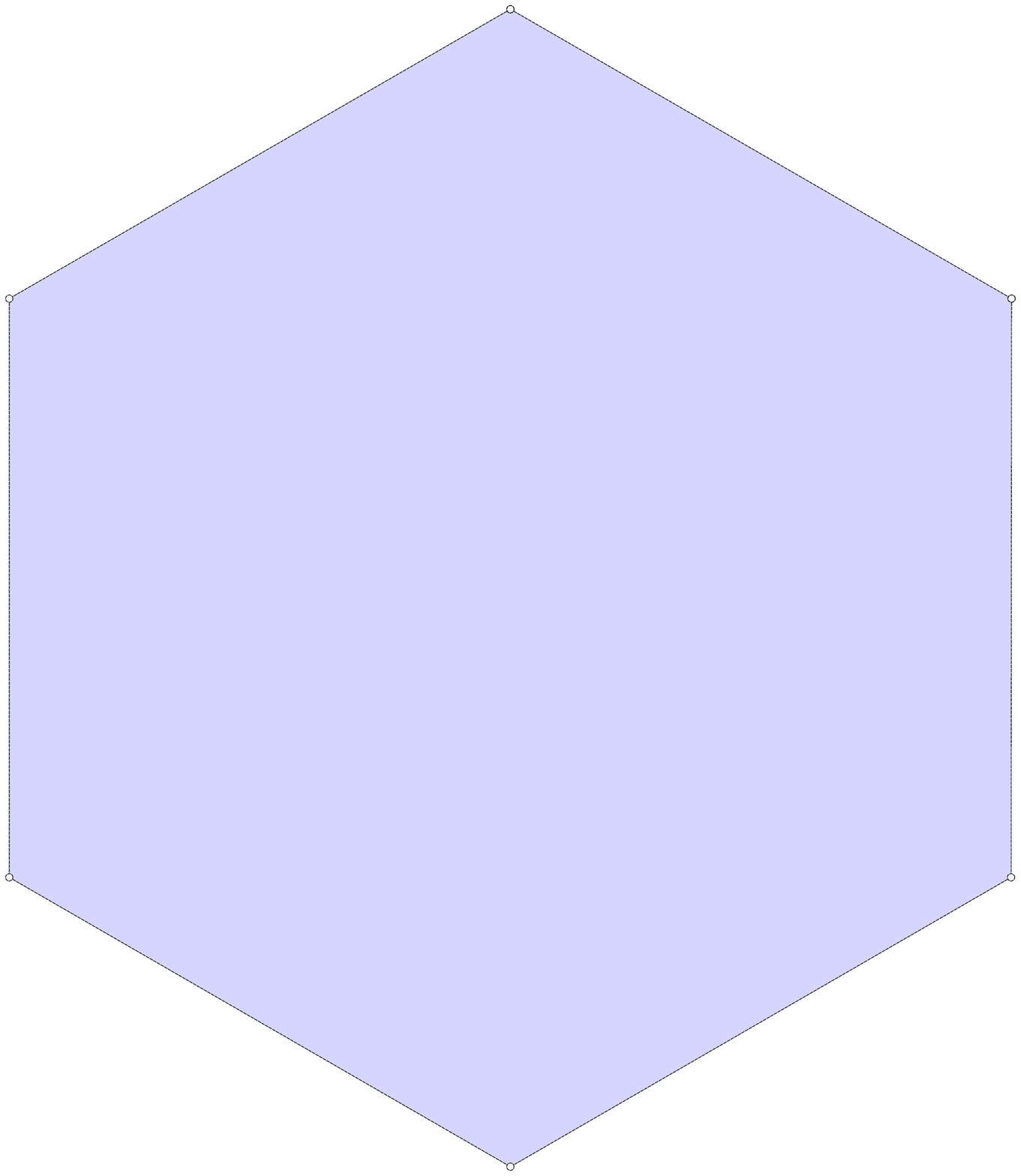,height=.5in,width=.4in}},\
\underset{(3^6)}{\epsfig{file=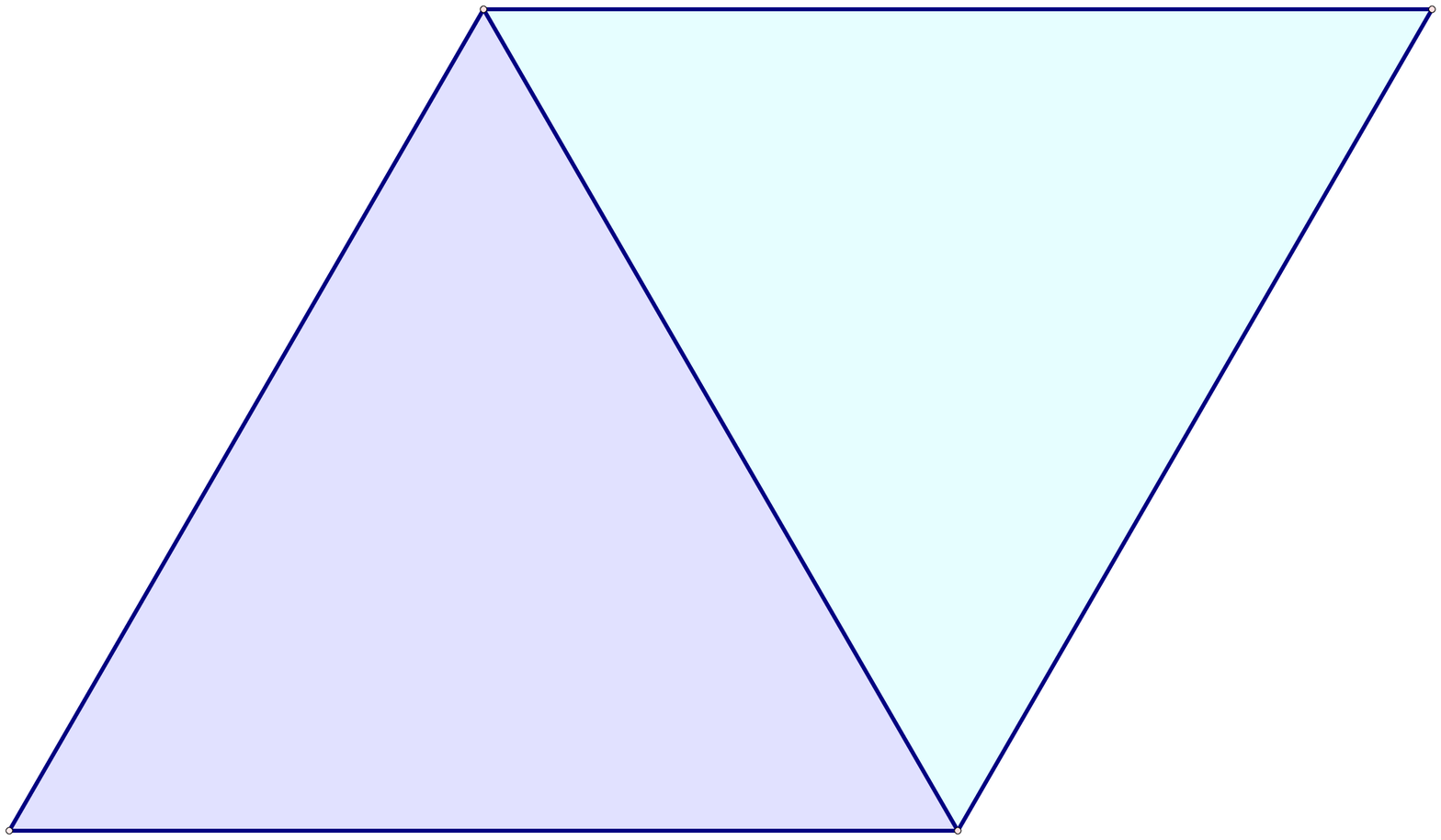,height=.3in,width=.5in}},\
\underset{(3^3,4^2)}{\epsfig{file=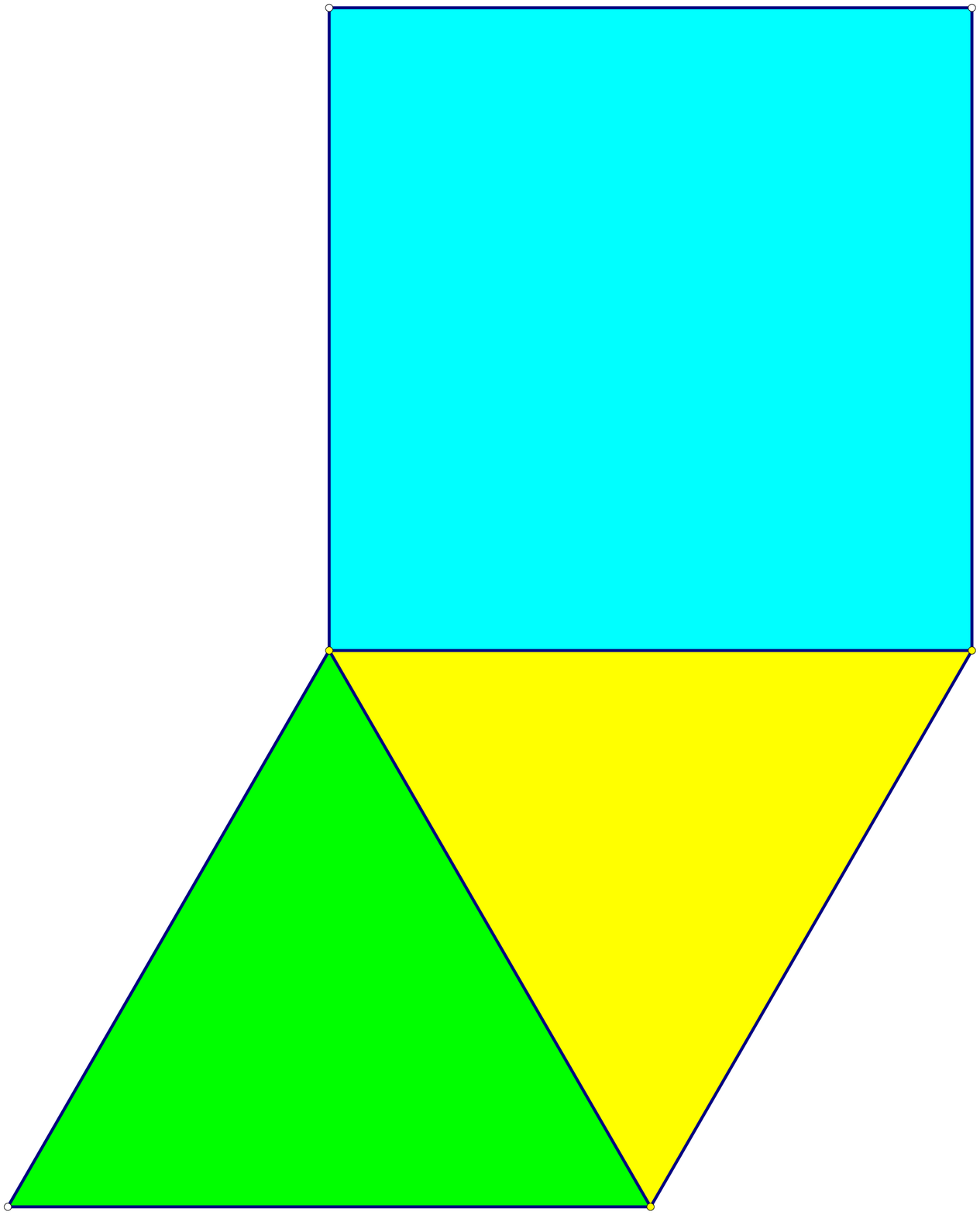,height=.5in,width=.4in}},\
\underset{(3,6,3,6)}{\epsfig{file=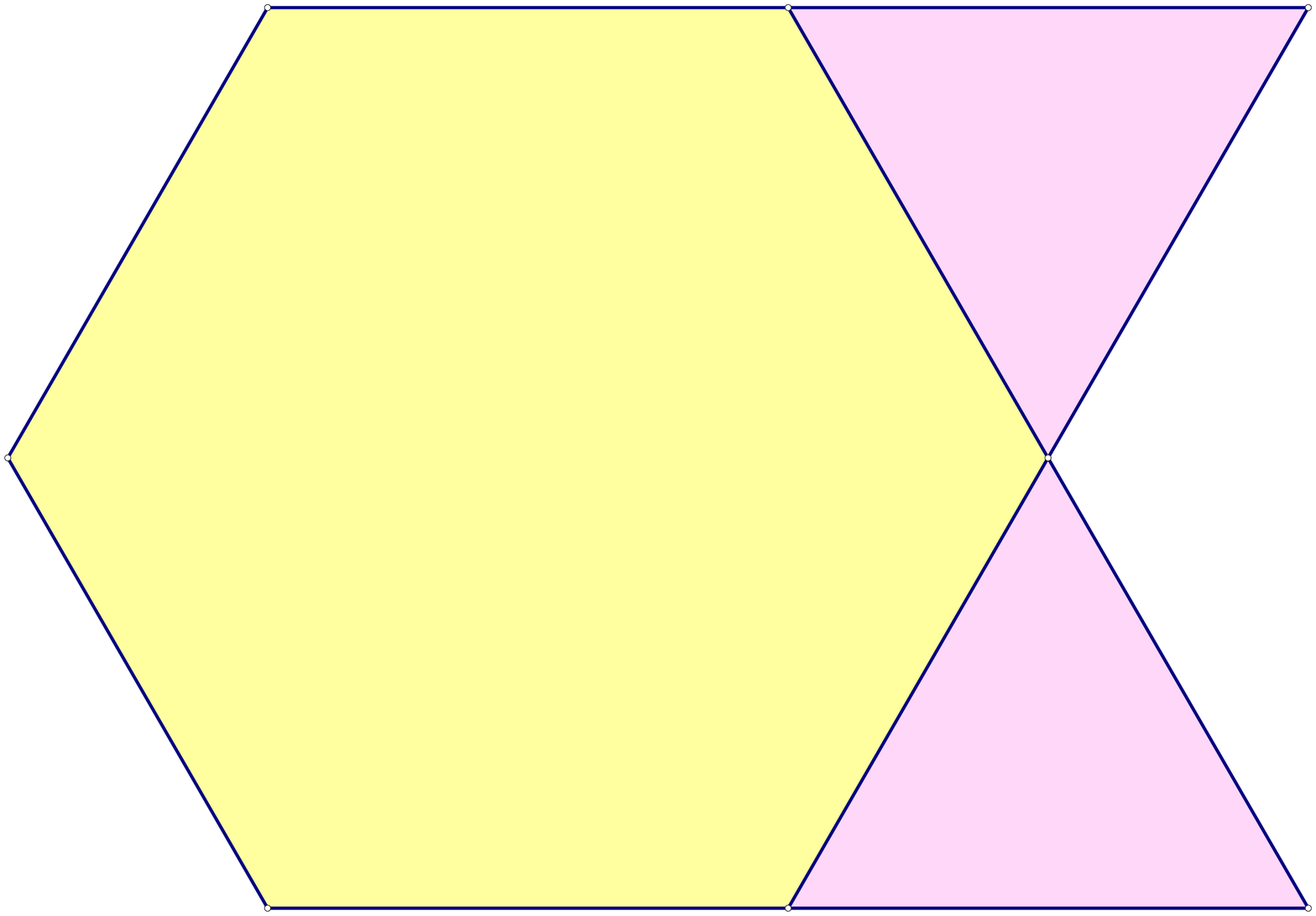,height=.4in,width=.5in}},
\underset{(3,4,6,4)}{\epsfig{file=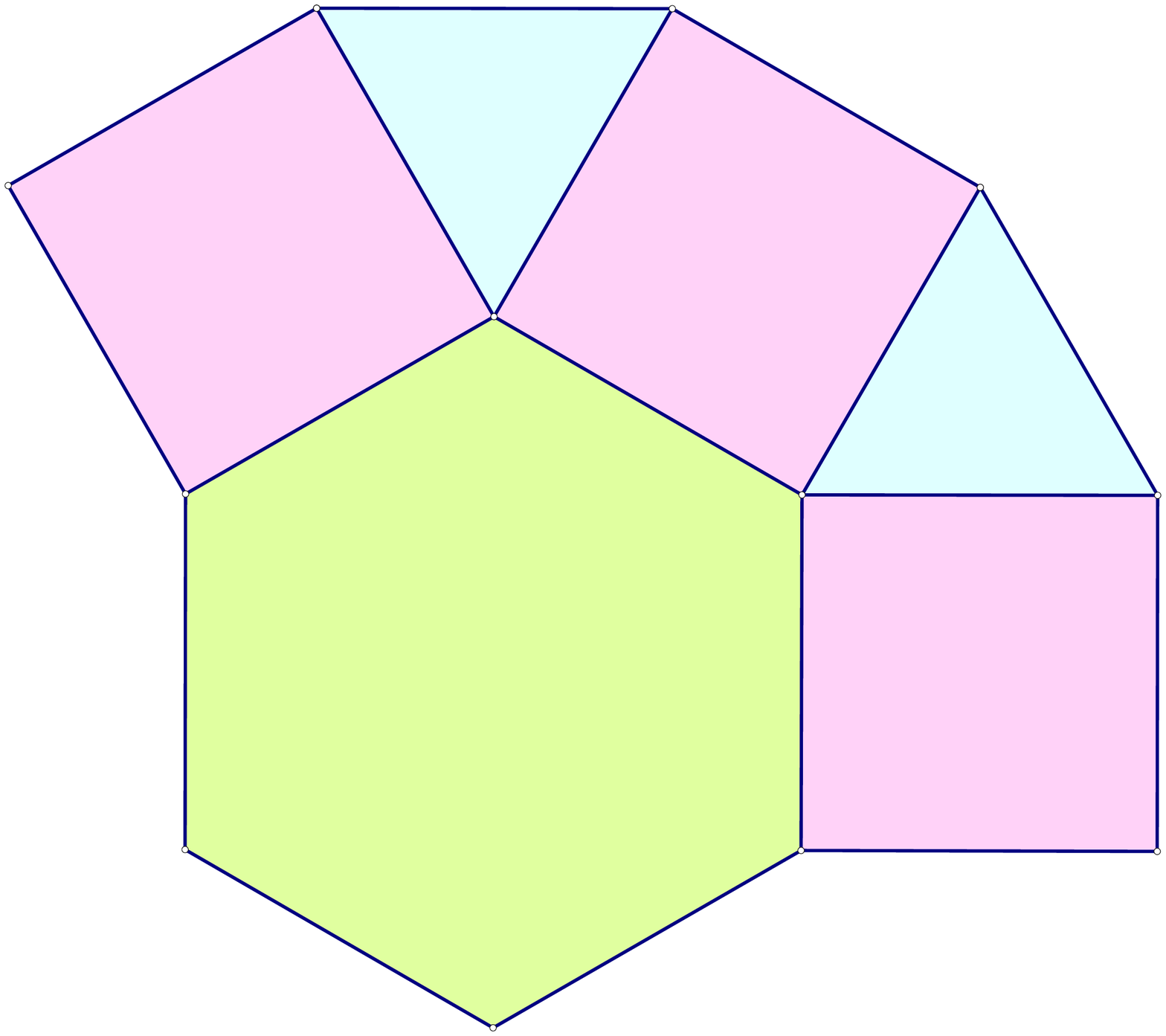,height=.5in,width=.5in}},\\ \\
\underset{(4,8^2)}{\epsfig{file=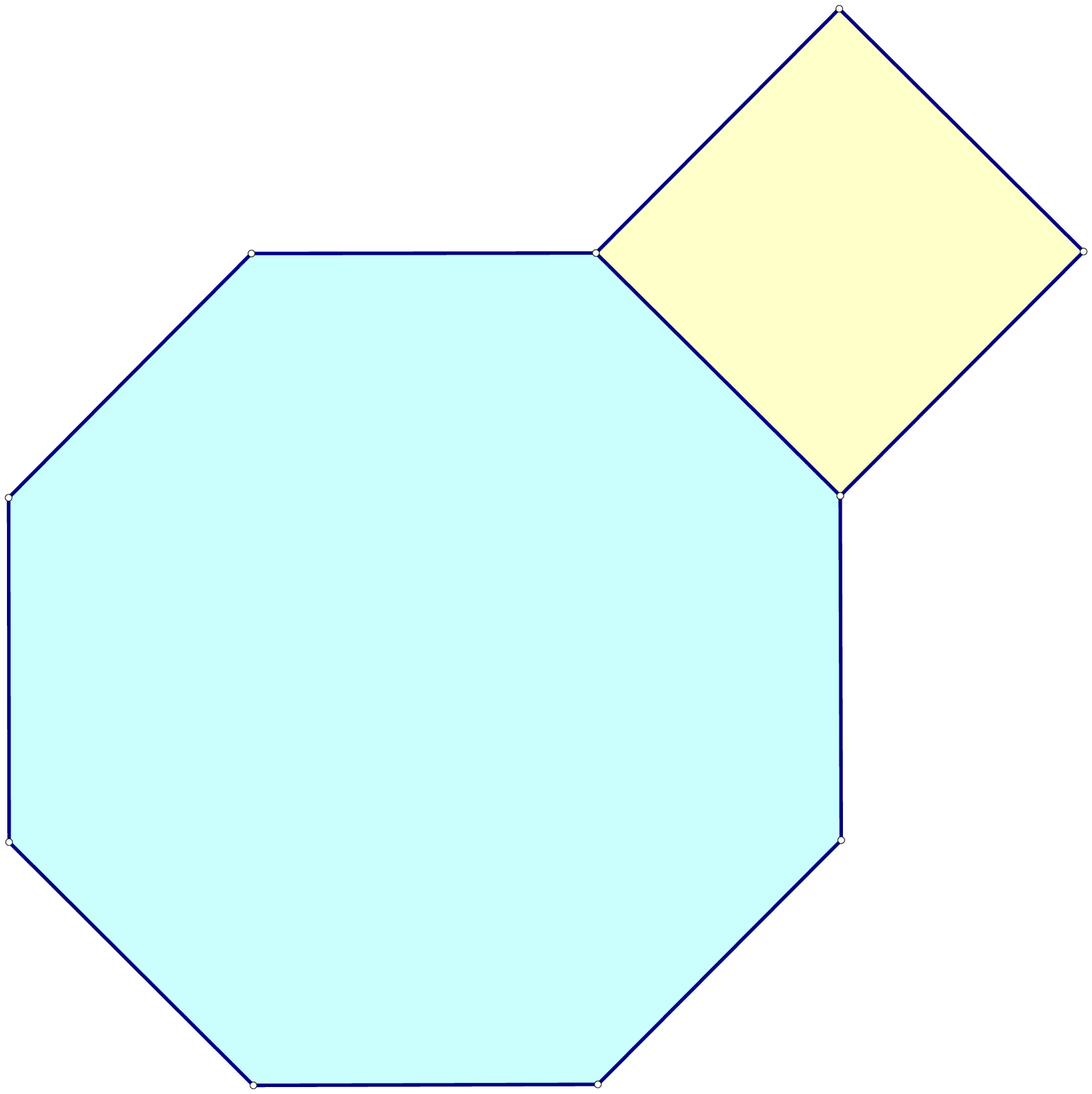,height=.5in,width=.5in}},
\underset{(4,6,12)}{\epsfig{file=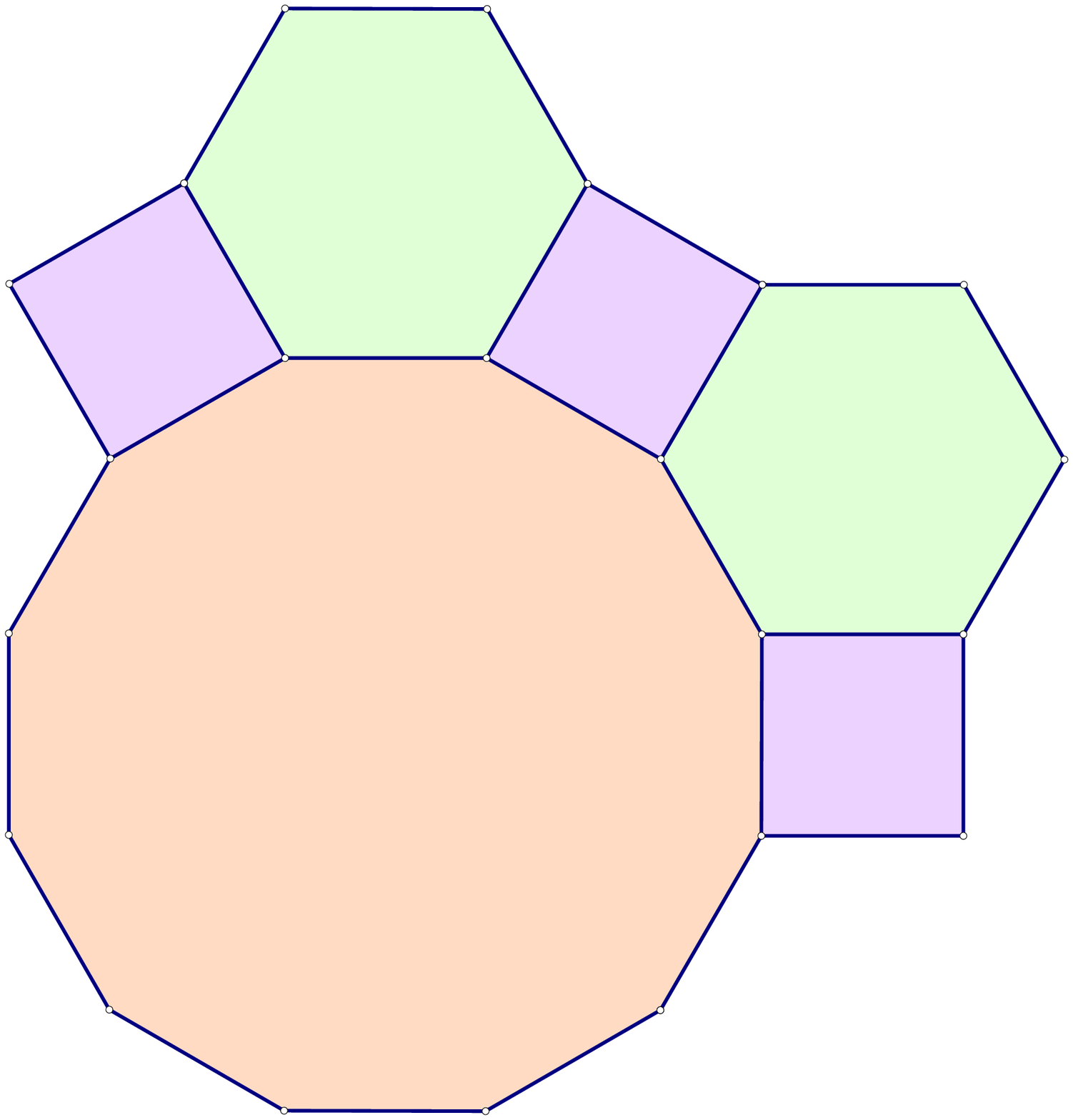,height=.5in,width=.5in}},\
\underset{(3,12^2)}{\epsfig{file=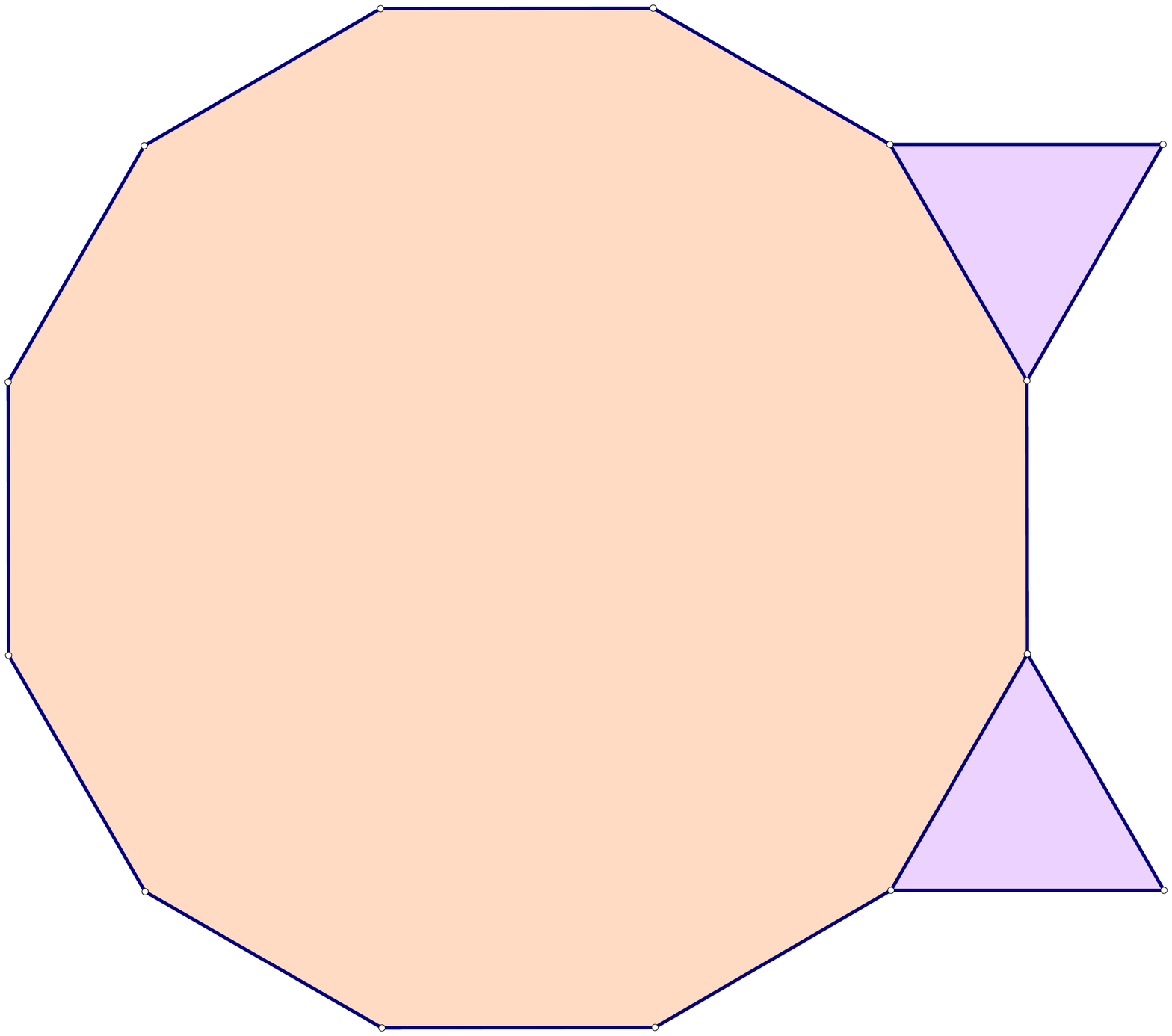,height=.5in,width=.6in}},\
\underset{(3^2,4,3,4)}{\epsfig{file=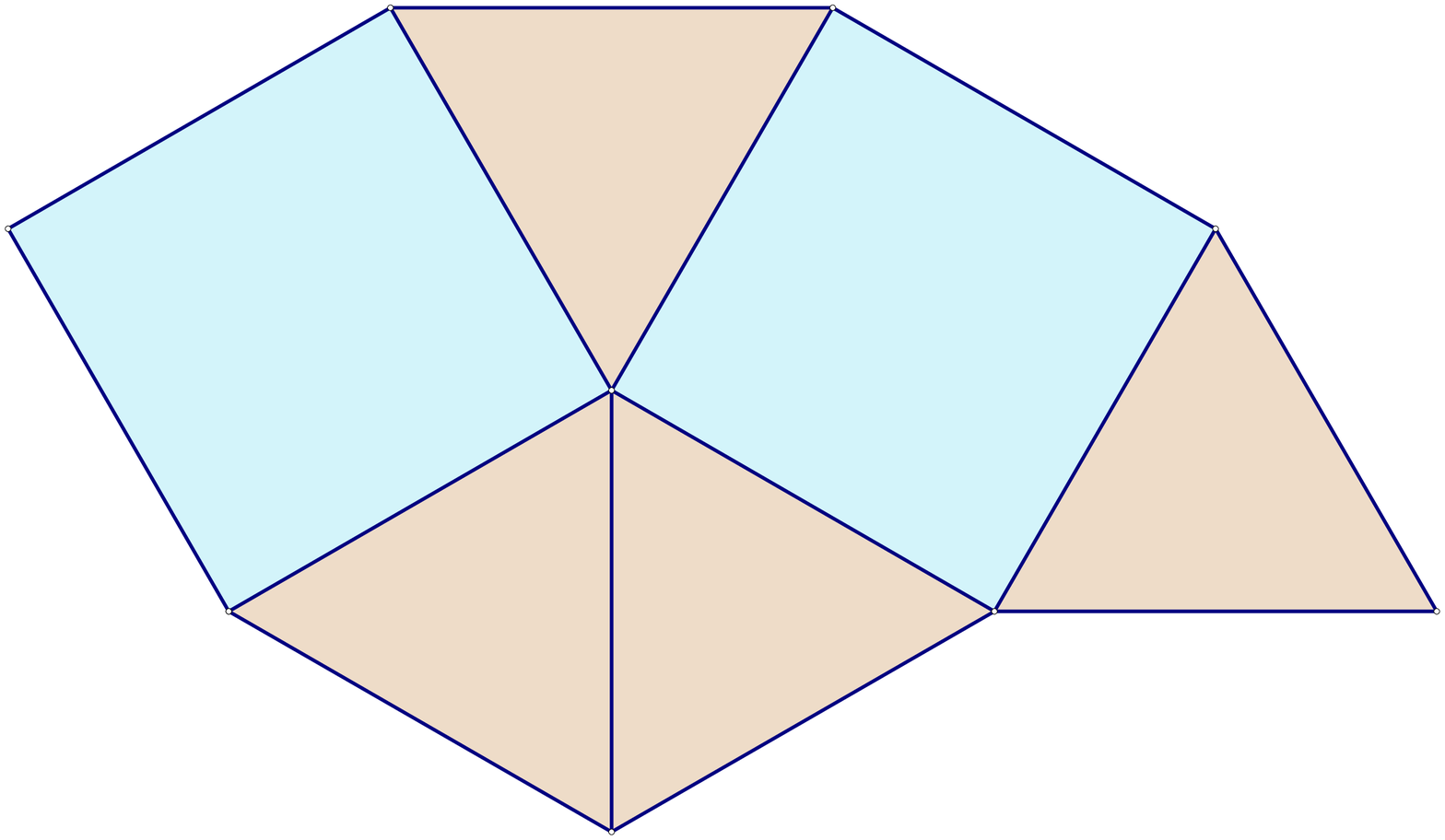,height=.5in,width=.8in}},\
\underset{(3^4,6)}{\epsfig{file=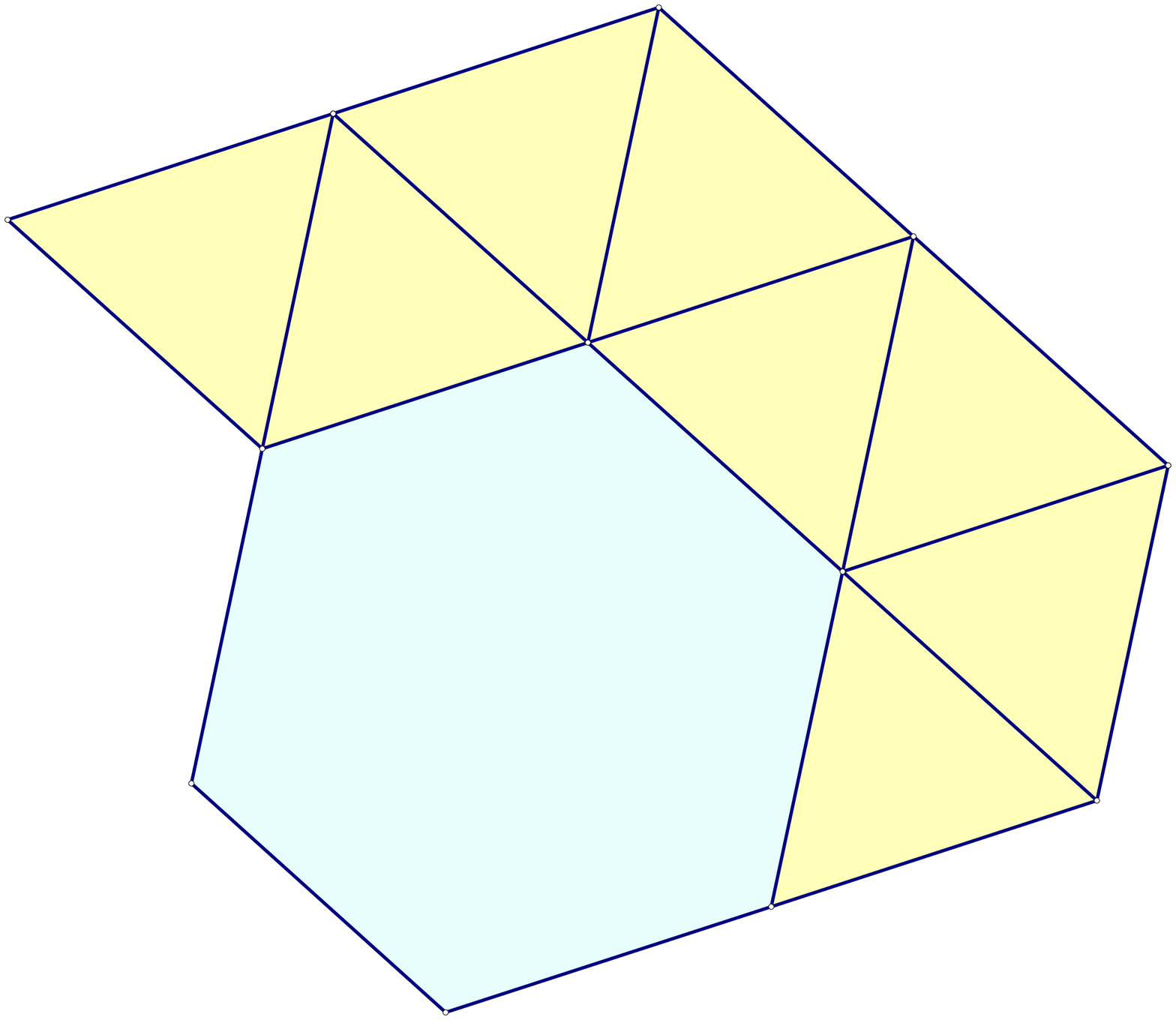,height=.5in,width=.6in}},\
\end{array}
\label{minimal}
$$\caption{Minimal clusters}\label{mintiles}
\end{figure}

By a {\it tiling} of the plane one
understands a countable union of closed sets (called tiles) whose
union is the whole plane and with the property that every two of
these sets have disjoint interiors. The term tessellation is a more
modern one that is used mostly for special tilings.  We are going
to be interested in the tilings in which the closed sets are either
all copies of one single regular convex  polygon (regular
tessellations) or several ones (semi-regular tessellations) and in
which each vertex has the same {\it vertex arrangement} (the number
and order of regular polygons meeting at a vertex). Also we are
considering the {\it edge-to-edge} restriction, meaning that every
two tiles either do not intersect or intersect along a common edge,
or at a common vertex. According to \cite{gs}, there are three
regular edge-to-edge tessellations and eight semi-regular
tessellations (see \cite{gs}, pages 58-59). The generic tiles in a
regular tessellation, or in a semi-regular one, are usually called {\it
prototiles}. For instance, in a regular tessellation the prototiles
are squares, equilateral triangles or regular hexagons. We will
refer to these tessellations by an abbreviation that stands for the
ordered tuple of positive integers that give the so called vertex
arrangement (i.e. the number of sides of the regular polygons around
a vertex starting with the smallest size and taking into account
counterclockwise order). The abbreviation is usually using the
convention with the powers similar to that used in the standard prime
factorization of natural numbers. So, the regular tessellation with squares
is referred to as $(4^4)$. We refer to Figure~\ref{mintiles} for the
rest of the notation.

Each such tessellation is {\it periodic} in the sense that there
exists a cluster of tiles formed by regular polygons, which by
translations generated by only two vectors, say $v_1$ and $v_2$,
covers the whole plane, and the resulting tiling is the given
tessellation (Figure~\ref{mintiles}).

\begin{figure}
$$\underset{(3^4,6)}{\epsfig{file=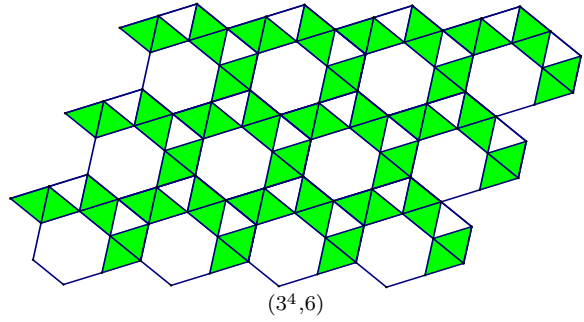,height=1.5in,width=3in}}$$
\caption{$G_{(3^4,6),4,3}$}\label{3^46mn84}
\end{figure}

For each such tessellation $\cal T$, we associate an infinite graph,
$G_{\cal T,\infty}$ in the following way. For each regular polygon
in the tessellation we have a vertex in $G_{\cal T,\infty}$, and
the edges in this graph are determined by every two polygons that share
a common edge within the tessellation. This way, for instance, for
the regular tessellation with squares $(4^4)$, we obtain what is
usually called the infinite {\it grid} graph. Also, if $m,n\in \Bbb N$ we
can construct a graph $G_{{\cal T},m,n}$ in the same way as before
using only $m$ copies of the cluster of tiles generating $\cal T$
shown in Figure~\ref{mintiles} in the direction of $v_1$ and $n$
copies in the direction of $v_2$. For ${\cal T}=(3^4,6)$, $m=4$ and
$n=3$ we obtain the graph generated by the polygons in the portion
of the tessellation $\cal T$ shown in Figure~\ref{3^46mn84}.

Clearly, the choice of the tiles on Figure~\ref{minimal} is not
unique but whatever one takes for these tiles it is not going to be
relevant in our calculation of densities. For each such graph as before, and one
of its vertices $v$, let us denote by $d(v)$ the number of adjacent
vertices to $v$, known as the degree of the vertex $v$. We define a
set of vertices $S$ to be {\it half-dependent} if for each vertex
$v\in S$ the number of adjacent vertices to $v$ that are in $S$ is
less than or equal to $\lfloor \frac{d(v)}{2}\rfloor$.

\begin{figure}[h]
\centering
$\underset{(a)}{\includegraphics[scale=1.2]{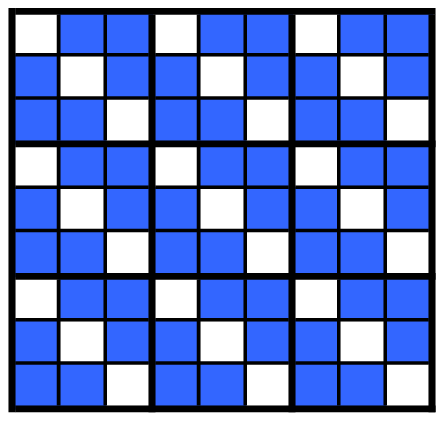}}$\ \ \
$\underset{(b)}{\includegraphics[scale=0.36]{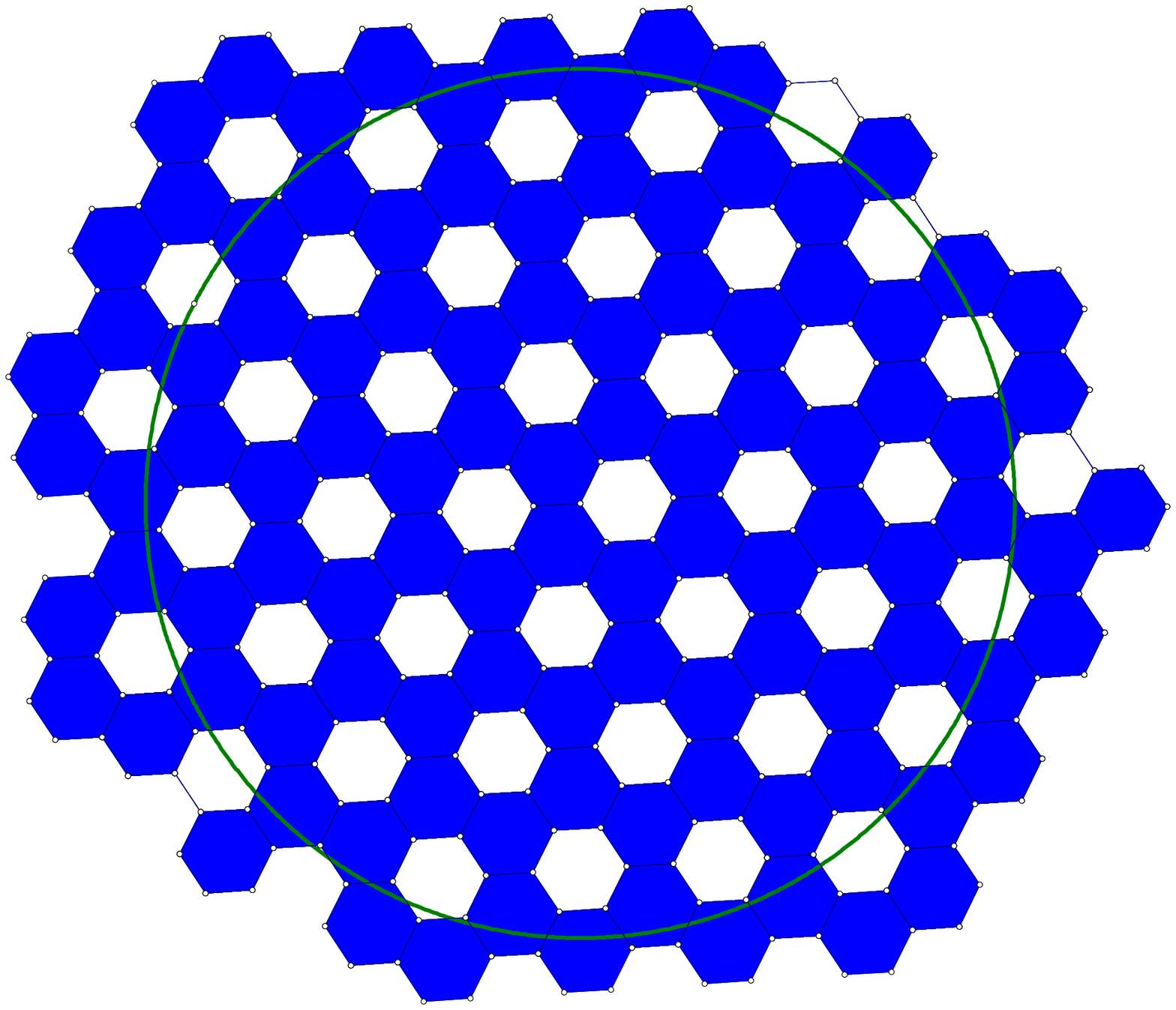}}$
\caption{Regular tessellations with squares and regular hexagons} \label{4444&666}
\end{figure}

Let $m,n\in \Bbb N$ be arbitrary, and for each graph $G_{{\cal
T},m,n}$ we denote by $\rho_{{\cal T}, m,n}$ the maximum cardinality of
a half-dependent set in $G_{{\cal T},m,n}$, divided by the number of
vertices of $G_{{\cal T},m,n}$. Hence one may consider the number

\begin{equation}\label{first}
\rho_{\cal T}=\limsup_{m,n \to \infty} \rho_{{\cal T}, m,n},
\end{equation}

\noindent which represents, heuristically speaking, the highest  proportion
in which one can distribute the vertices in a half-domination set.
 For instance we will show that in the case of a regular
tessellation with regular hexagons (see Figure~\ref{4444&666})(b)
the number defined above is $\frac{2}{3}$.

\noindent In the  regular tessellation with squares we have shown in
\cite{hiw} that the number defined in (\ref{first}) is also
$\frac{2}{3}$ (Figure~\ref{4444&666}(a)).
In this paper we are interested in the values of $\rho_{\cal T}$ for these types of tessellations.

\section{Various techniques}
First we are going to use one of main ideas in \cite{ipw} and \cite{hiw}, and some classical linear optimization techniques.
\subsection{Integer Linear Programming Systems}

The regular tessellation $(6^3)$ can be treated simply in the following way. We are going to work with the graph
 $G_{6^3,n,n}$ obtained from translating $n$-times a prototile so that the each two neighboring  tiles share an edge and then translate the whole
 row of n-tiles in such a way that two neighboring rows fit together to give the tessellation of a rhombic n-by-n region. We assign to each tile
 a variable $x_{i,j}$ which can be $0$ or $1$: if a tile is part of the half-dominating set (colored blue in Figure~\ref{4444&666}b) of maximum cardinality then
 its variable $x_{i,j}$ is equal to $1$, and if it is not,  $x_{i,j}=0$.  For most of the vertices   of $G_{6^3,n,n}$ the degree is $6$ so each vertex in the half-dominating set, denoted by $V$,  must have at most $3$ other vertices adjacent to it which are in $V$. Let us denote by $N(i,j)$ the adjacent indices to the vertex indexed with (i,j). We can write the half-domination condition as

 $$x_{i,j}^{\star}:=\underset{(k,l)\in N(i,j)}{ \sum }  x_{k,l}\le 3\ \text{ for all (i,j) with } 1<i,j<n, (i,j)\in V.$$

\noindent  For a vertex not in $V$ we simply have no restriction on the above sum (maximum is $6$). The trick is to write an inequality that encompasses both situations. In this case,
 we see that the following inequality accomplishes exactly that:

 $$3x_{i,j}+\underset{(k,l)\in N(i,j)}{ \sum }  x_{k,l}\le 6\ \text{ for all (i,j) with } 1<i,j<n.$$

 \noindent For the vertices on the boundary we have similar inequalities.
 Summing all these inequalities up, we obtain, $3|V'|+6|V'|+ \Sigma\le  6n^2$, where $V'$ are the vertices in $V$ in the interior and
 $\Sigma$ are the number of those on the boundary. Since clearly $\lim \frac{\Sigma}{n^2}=0$ we see that $\lim\frac{|V|}{n^2}\le \frac{2}{3}$.
 In Figure~\ref{4444&666}b, we exemplify an arrangement which accomplishes this density. Hence, $\rho_{6^3}=\frac{2}{3}$.

\vspace{0.1in}

In the case of the tessellation with equilateral triangles, the inequality above changes into

 $$2x_{\triangle}+\underset{\triangle' \in N(\triangle)}{ \sum }  x_{\triangle'}\le 3\ \text{ for all\ } \triangle s \ \text{ not at the boundary}.$$

\noindent The argument above gives the inequality $\rho_{3^6}\le \frac{3}{5}$. The arrangement in Figure~\ref{333333} has a density which is
 only $\frac{9}{16}$. In this case it is not possible to achieve the density $\frac{3}{5}$.

\subsection{Toroidal type graphs}

One way to avoid to deal with the boundary tiles, is to form  toroidal type graphs obtained from tessellations.
We have shown in \cite{ipw} that this does not affect the maximum density, in the sense that both situations tend at the limit to the same density value.
Let us begin with the case ${\cal T} =3^6$. For $n\in \mathbb N$, we take the parallelogram which gives the graph $G_{3^6,n,n}$ and identify the
opposite edges (without changing the direction). This gives rise to a similar graph which we will denote by $T_{3^6,n}$. This graph is regular and the equations
 we get can be easily described and implemented in LPSolve. Let us introduce the variables $x_{i,j,k}$, $i,j\in \{0,1,2,...,n-1\}$ and $k\in \{1,2\}$, in the following way. The index $(i,j)$ refers to the translation of the minimal  parallelogram (Figure~\ref{mintiles})  $i$ places in the horizontal direction and $j$
places in the $60^{\circ}$ direction. The index $k$ corresponds to the lower ($k=1$) and upper equilateral triangle within this parallelogram ($k=0$). We get a number of $2n^2$ vertices of this regular graph. The optimization equations are simply

\begin{equation}\label{equations333333}
\begin{cases}
2x_{i,j,1}+x_{i,j,2}+x_{u,j,2}+x_{i,v,2}\le 3, \ where \ u\equiv i-1\ (mod\ n), \\ \\ \ \ \  v\equiv j-1\ (mod \ n), u,v\in \{0,1,2,...,n-1\}\\ \\
2x_{i,j,2}+x_{i,j,1}+x_{u,j,1}+x_{i,v,1}\le 3, \ where \ w\equiv i+1\ (mod\ n), \\ \\ \ \ \ t\equiv j+1\ (mod \ n), w,t\in \{0,1,2,...,n-1\}, i,j \in \{0,1,2,...,n-1\}.
\end{cases}
\end{equation}

\noindent We denote by $\rho_{3^6,n}$ the best densities for this graph.   Observe that  $\rho_{3^6,n}\le \rho_{3^6,m}$ if $n$ divides $m$.
Also, as proved in \cite{ipw}, we have $\rho_{3^6}=\underset{n} {\sup} \rho_{3^6,n}$. We have the following densities:

\centerline{
\begin{tabular}{|c|c|c|c|c|c|c|c|c|c|}
  \hline
  n & 1 & 2 & 3 & 4 & 5 & 6& 7 & 8 & 9 \\
  \hline
  $\rho_{3^6,n}$ & 0 & $\frac{1}{2}$ & $\frac{5}{9}$ & $\frac{9}{16}$ & $\frac{14}{25}$ & $\frac{5}{9}$ &  $\ge \frac{27}{49}$& $\ge \frac{9}{16}$ & $\ge \frac{5}{9}$ \\
  \hline
\end{tabular}}

\vspace{0.1in}

\begin{figure}[h]
\centering
$\underset{Density \ \frac{9}{16}}{\includegraphics[scale=0.3]{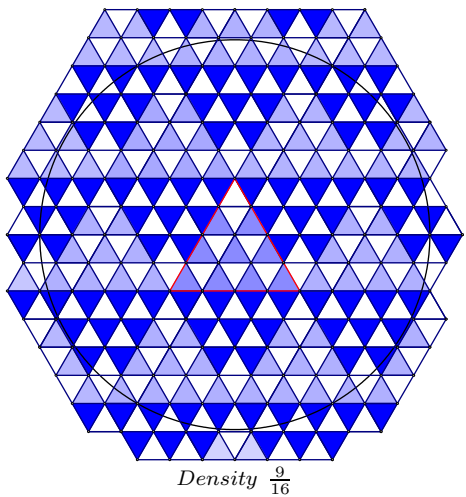}}$
\caption{Regular tessellation with equilateral
triangles} \label{333333}
\end{figure}

\begin{figure}[h]
\centering
\includegraphics[scale=0.36]{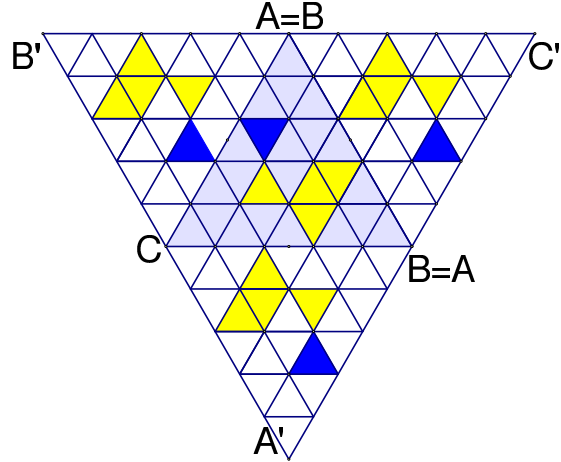}
\caption{Klein type toroidal identification} \label{kleintoroidal}
\end{figure}

\noindent Unfortunately, LPSolve takes too long to finish the analysis for $n\ge 7$ (about 100 variables). In order to reduce the number of variables,
and of course taking advantage of the rotational symmetry, one can try to do a different matching of the boundaries as in Figure~\ref{kleintoroidal}.
This corresponds to rotations of $180^{\circ}$ around the midpoints of the sides
for the basic triangle of tiles $ABC$. For instance, the corner tiles are all adjacent to each other.

The equations that we need to use, are those given by (\ref{equations333333}) and in addition

\begin{equation}\label{kleinequations}
\begin{array}{c}
x_{i,j,1}=x_{u,v,2}\ where \ u\equiv n-1-i\ (mod\ n), \  and\ \    v\equiv n-1-j\ (mod \ n),\\ \\  u,v\in \{0,1,2,...,n-1\}, \text{for all}\  i,j\in \{0,1,2,...,n-1\}.
\end{array}
\end{equation}
\noindent We will refer to the graph obtained by these identifications as $K_{3^6,n}$. This graph is a regular graph with $n^2$ vertices and degree 3.
 All the densities are the same as in the table above but LPSolve finishes in a reasonable time and the arrangement in Figure~\ref{333333} satisfies this
new restriction. In the case $n=4$ we can see that the arrangement in Figure~\ref{333333} is also the best because if we could add another triangle, that would increase the density to $\frac{5}{8}>\frac{3}{5}$, which we know is not possible by the upper bound established earlier. It is important that the best density in the case $K_{3^6,8}$
is the same as in the case $K_{3^6,4}$. We conjecture that the best density is given by such a matching of two toroidal graphs one having the dimension double the dimension of the other.

\subsection{Upper bounds}

For the problem above we can adopt the method of weighted objective function as described in \cite{hiw}. We used the
weights all equal to $1$ in the interior and zero on the boundary for $K_{3^6,13}$. The upper bound obtained is $\rho_{3^6}\le \frac{70}{121}\approx 0.59375$.
This is within $\frac{31}{1936}\approx 0.0160124$ to $\frac{9}{16}$ which makes it plausible that  $\rho_{3^6}=\frac{9}{16}$.

For the graph $K_{3^6,n}$, we may add to the system of inequalities  (\ref{equations333333}) and (\ref{kleinequations}) the conditions
$0\le x_{i,j,1}, x_{i,j,2}\le 1$ and $\sum_{i,j} (x_{i,j,1}+x_{i,j,2})=k$. These inequalities describe a polytope in $n^2$ dimensions. Finding the maximum cardinality of a half dependent set is equivalent to finding the smallest $k$ for which there is no lattice point in the corresponding polytope. There exists a theory which counts the number of lattice points in polytopes which was started by Eugene Ehrhart in 1960's (see \cite{abarvinok} and \cite{beckAndRobins2007textbook}). The theory simplifies significantly if the polytope has integer vertices. Unfortunately our polytope has
rational vertices.  Theoretically, there exists a quasi-polynomial $P(t,k)$ of degree $n^2$ which counts the lattice points contained in the dilation of the polytope by $t$.
We want the smallest $k_n$ such that $P(1,k_n+1)=0$. Let us make the observation that  $P(1,k_n)$ should by relatively a big since the system is invariant under translations (mod $n$) and under rotations. So, one may expect $P(1,k_n-1)\approx 3n^2$. There exists several programs which calculate this polynomial from the coefficients of the system of inequalities which define the polytope. One of these programs is called LattE and it written by Matthias  K$\rm \ddot{o}$ppe and Jes$\rm \acute{u}$s A. De Loera. This method remains to be implemented and investigated in a different project.  Also, it seems like the number of variables that can be used in this program is also by about 100. However, we think that the method fits very well with the toroidal formulation. If the Ehrhart polynomial can be computed in dimensions $k$ and $2k$, it may give more information of how to find the best density.
We will be using other methods which are discussed in the next sections.

\section{Semi-regular tessellations}

Working with the graphs generated by semi-regular tessellations is more challenging when it comes to implementing
the problem into LPSolve. There are also advantages here since the systems have less symmetry and somehow this is a plus for
the optimization programs to arrive faster at a maximum.

\subsection{The case ${\cal T}=(3^3,4^2)$}

\begin{figure}[h]
\centering
$\underset{(a)\ Density \ \frac{7}{12}}{\includegraphics[scale=0.36]{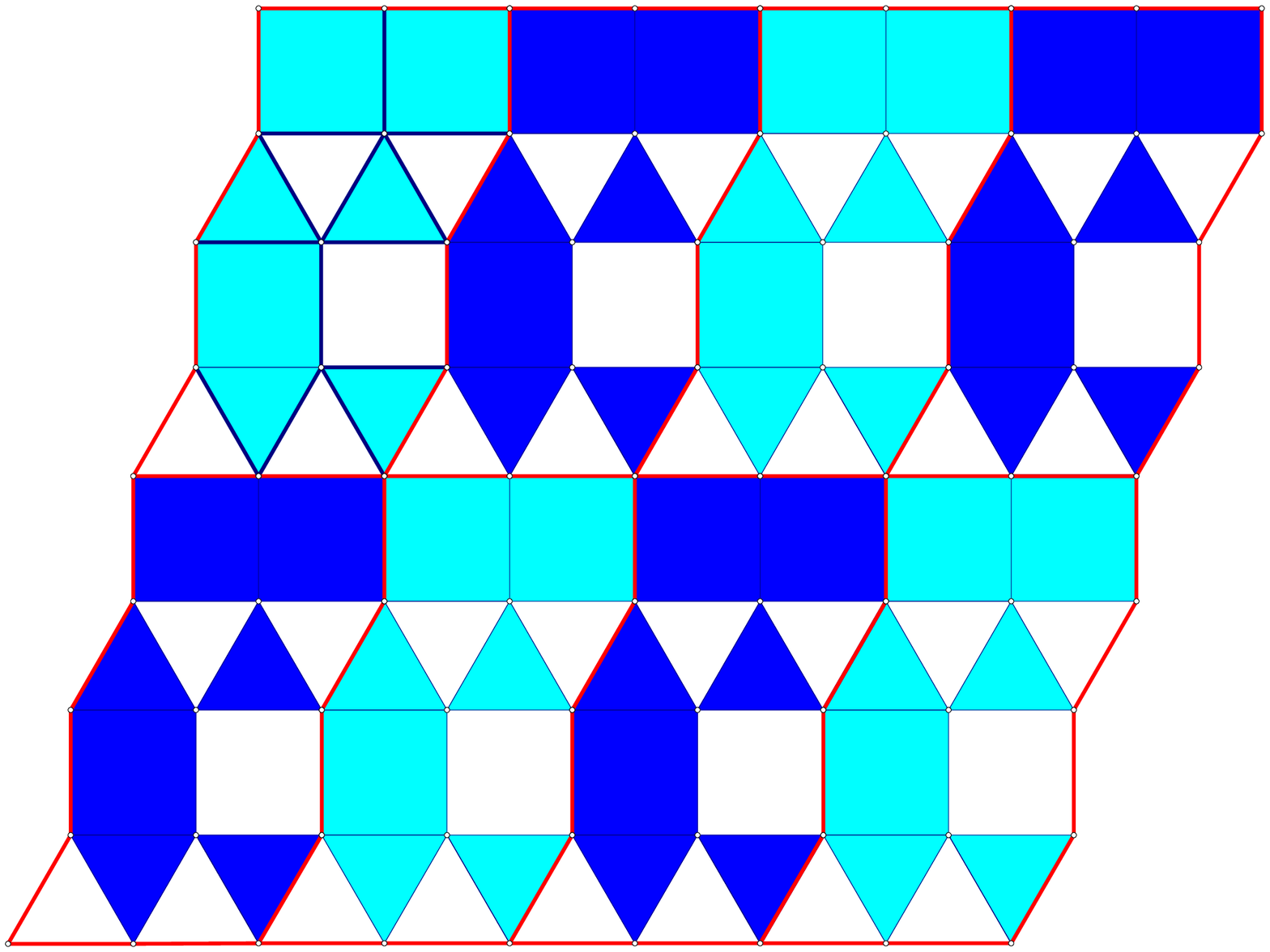}}\
\underset{(b)\ Density \ \frac{11}{18}}{\includegraphics[scale=0.3]{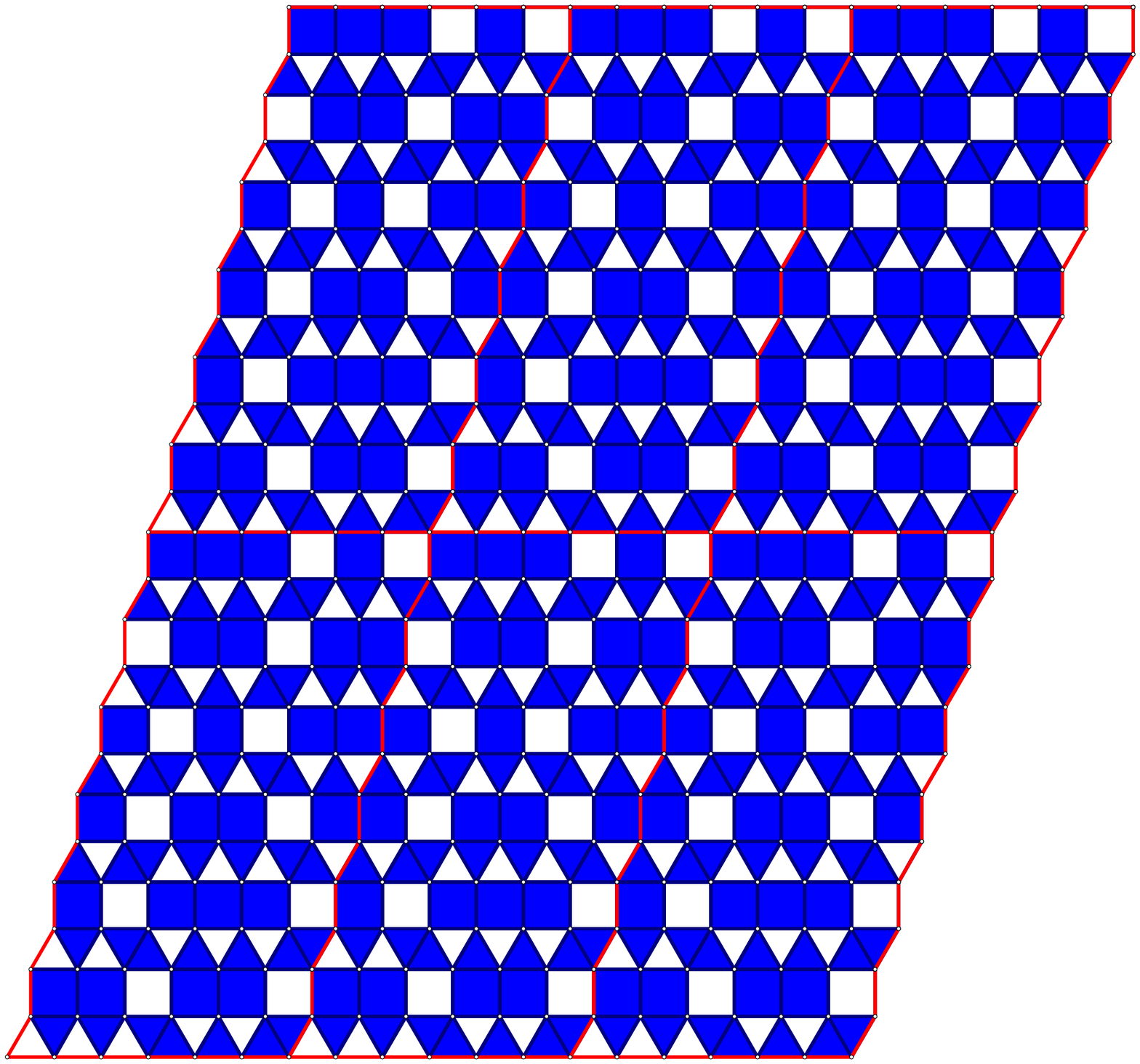}}$
\caption{Some arrangements for ${\cal T}=(3^3,4^2)$ } \label{arrangementin33344}
\end{figure}

\n In Figure~\ref{arrangementin33344} (a) we see an arrangement of a half-domination set with a density of $\frac{7}{12}$, which is only $\frac{1}{28}$
off of the upper bound in the next theorem.

\begin{theorem}\label{bestupperbound33344}
The half-domination density for the tessellation ${\cal T}=(3^3,4^2)$ satisfies
$$ \rho_{(3^3,4^2)}\le \frac{13}{21}.$$
\end{theorem}

\n \proof. Let us consider a toroidal graph induced by $G_{{\cal T},m,n}$ with $m$ and $n$ big and a half domination set corresponding to it.
We observe that there are $\frac{3mn}{3}=mn$ squares and $\frac{2}{3}(3mn)=2mn$ equilateral triangles in $G_{{\cal T},m,n}$.
As usual let us consider the variables $x_v$ for each vertex $v$ in this graph defined to be 1 or 0 as being in the domination set or not.
Also, we denote by $x_v^{*}$ the sum of the variables $x_w$ corresponding to the adjacent vertices $w$ of $v$. For a vertex $v$ corresponding to
a square we have (see Figure~\ref{arrangementin33344})

\begin{equation}\label{firsteq}
2x_v+x_v^{*}\le 4.
\end{equation}

\n We denote by $T$ the sum of all $x_v$ over all vertices corresponding to triangles and by $S$ the sum of all $x_v$ over all vertices corresponding to squares.
If we sum up all equalities (\ref{firsteq}) over all the squares we get:

$$2S+\underset{\text {v for a squares}}{\sum x_v^{*}} \le 4mn \ \Rightarrow \ 2S+(2S+T)\le 4mn\ \text{or}\  4S+T\le 4mn.$$

\n For a vertex $v$ corresponding to
an arbitrary triangle we have

\begin{equation}\label{secondteq}
2x_v+x_v^{*}\le 3,
\end{equation}

\n which gives, as before, if summed up over all of the triangles:

$$2T+\underset{\text {v for a trianlge}}{\sum x_v^{*}} =6mn\ \Rightarrow \ 2T+(2S+2T)\le 6mn\ \text{or}\ 2S+4T\le 6mn.$$

\n If we let $x=\frac{S}{3mn}$ and $y=\frac{T}{3mn}$, we need to maximize $x+y$ and, as we have seen above, $x$ and $y$ are subjected to the restrictions

$$\begin{cases} 4x+y\le \frac{4}{3}\\ \\
2x+4y\le 2.
\end{cases}$$

\begin{figure}[h]
\centering
\includegraphics[scale=0.36]{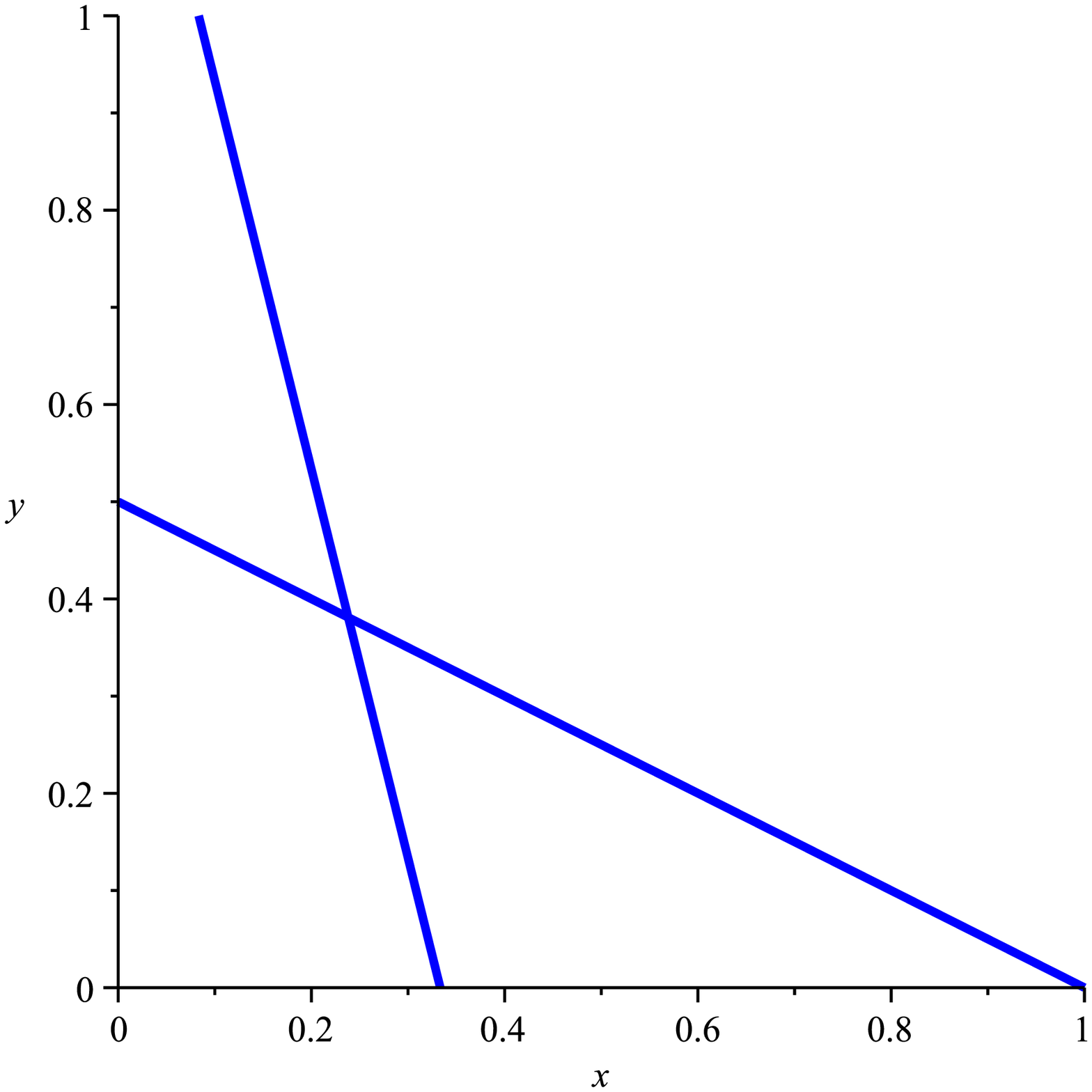}\ \ \  \includegraphics[scale=0.36]{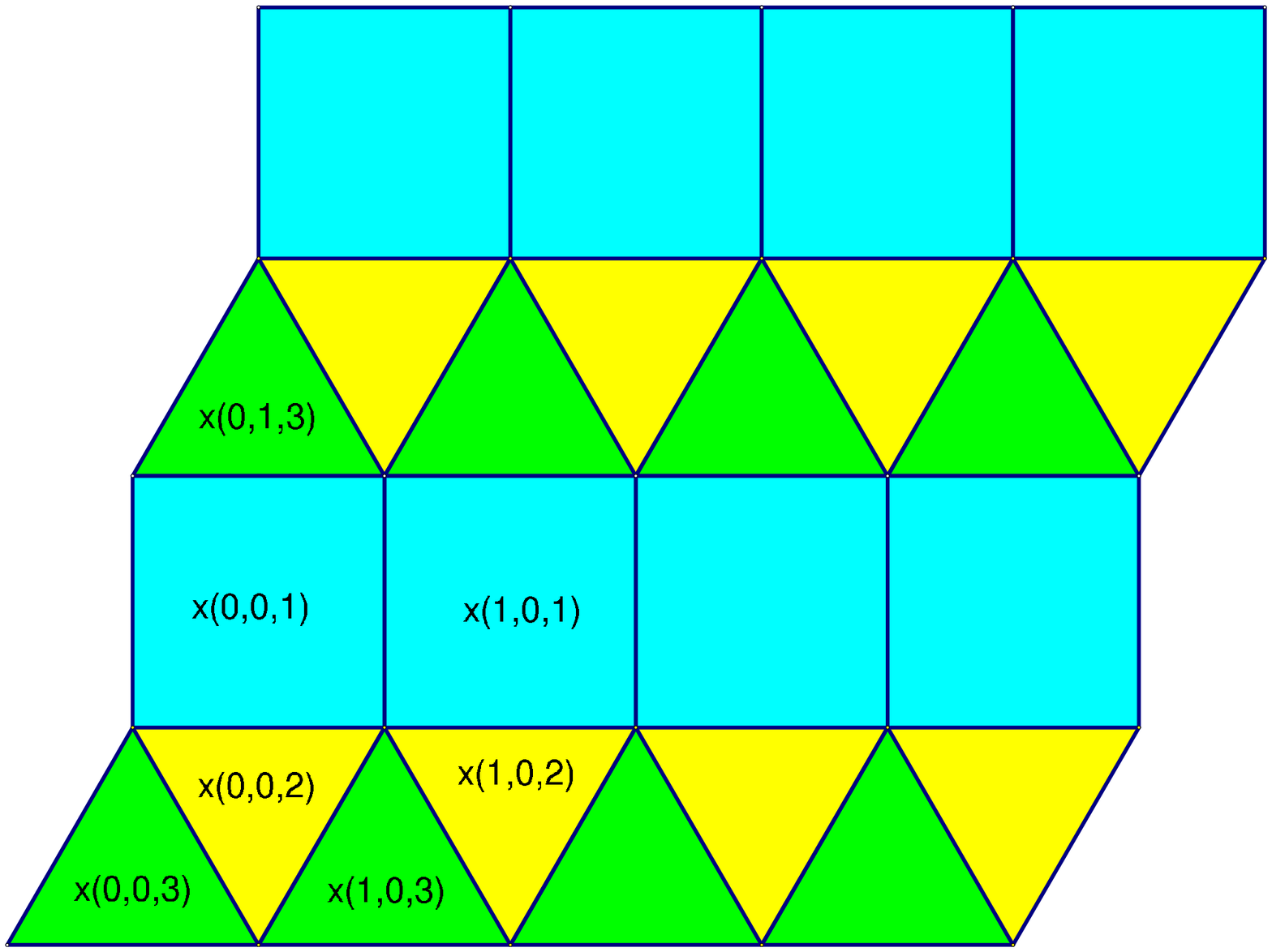}
\caption{The two lines and the variables convention } \label{twolines}
\end{figure}

\n The usual maximization technique (see Figure~\ref{twolines}) gives that for real numbers $x$ and $y$, the maximum of
$x+y$ is attained for $x=\frac{5}{21}$ and $y=\frac{8}{21}$ with $x+y=\frac{13}{21}$. This proves that
$$ \rho_{(3^3,4^2)}=\limsup_{m,n\to \infty} \frac{S(m,n)+T(m,n)}{3mn}\le \frac{13}{21}.\ \ \ \bsq$$

\subsection{Toroidal examples for ${\cal T}=(3^3,4^2)$}

The number of variables grows rapidly with $m=n$. For $n=7$ we have $3n^2=147$ variables and this seems to be a good bound for what one can obtain with
LPSolve. We use variables $x_{i,j,k}$ (having values 0 or 1) with $i=0,1,2,...,n-1$, $j=0,1,2,...,m-1$ and $k=1,2,3$ with $x_{i,j,1}$ for a square, $x_{i,j,2}$  and $x_{i,j,3}$ for
the triangles in the cluster of minimal tiles as in the Figure~\ref{twolines}. The  system of inequalities can be written as follows

$$\begin{cases}
2x_{i,j,1}+x_{i,j,2}+x_{u,j,2}+x_{v,j,1}+x_{i,t,3}\le 4, \ where \ u\equiv i-1\ (mod\ n),\  v\equiv i+1\ (mod \ n), \\ \\
w\equiv j-1\ (mod\ m),\  t\equiv j+1\ (mod \ m), \ u,v\in \{0,1,2,...,n-1\},w,t\in \{0,1,2,...,m-1\}\\ \\
2x_{i,j,2}+x_{i,j,1}+x_{u,j,3}+x_{v,j,3}\le 3,\ \text{and}\ 2x_{i,j,3}+x_{i,j,2}+x_{u,j,2}+x_{i,w,1}\le 3.
\end{cases}$$

\n The best densities we have gotten are listed next:

\vspace{0.1in}
\centerline{
\begin{tabular}{|c|c|c|c|c|c|c|c|}
  \hline
  n & 1 & 2 & 3 & 4 & 5 & 6& 7  \\
  \hline
  $\rho_{3^34^2,n}$ & 0 & $\frac{1}{2}$ & $\frac{5}{9}$ & $\frac{7}{12}$ & $\frac{3}{5}$ & $\frac{11}{18}$ &  $\ge \frac{4}{7}$ \\
  \hline
\end{tabular}}
\vspace{0.1in}

\n Of these arrangements, the case $m=n=6$ (Figure~\ref{arrangementin33344} (b)) gives the highest  density which is within $\frac{1}{126}$ to the upper bound shown above.

The same method used before for finding upper bounds does not give a better  bound as the one we have proved in Theorem~\ref{bestupperbound33344}.
For $m=n=7$ we let $x_{i,0,2}=x_{i,0,3}=0$ and $x_{0,j,1}=x_{0,j,2}=x_{0,j,3}=0$ with $i \in \{0,1,2,...,n-1\}$, $j\in \{0,1,2,...,m-1\}$ and LPSolve
gives a maximum of $72$. So, the upper bound is $\frac{72}{147-(7(2)+7(3)-2)}=\frac{12}{19}>\frac{13}{21}$.

In \cite{hiw} we have introduced the concept of deficiency function, $\delta_{i,j,k}$, and global deficiency $\Delta$ of an arrangement. Let us see how this works in this situation.
We define

$$\delta_{i,j,k}=
\begin{cases}\ \text{if} \ x_{i,j,k}=1\ \begin{cases}  x^{*}_{i,j,1}-2\ \text{if}\ k=1\\ x^{*}_{i,j,k}-1\ \text{if}\ k=1\ \text{or}\ 2,
\end{cases}\\ \\
\text{if}\  x_{i,j,k}=0\ \begin{cases}  x^{*}_{i,j,1}-4\ \text{if}\ k=1\\ x^{*}_{i,j,k}-3\ \text{if}\ k=1\ \text{or}\ 2,
\end{cases}
\end{cases}\ \text{and}$$

$$ \Delta=\frac{1}{|V|}\underset{(i,j,k)\in V}{\sum} \delta_{i,j,k}, \text{where} \ V \ \text{is the set of vertices.}
$$

\n We observe that the arrangement in Figure~\ref{arrangementin33344} (a) has $\Delta=\frac{-2}{12}=-\frac{1}{6}$ and the arrangement in
Figure~\ref{arrangementin33344} (b) has $\Delta=-\frac{1}{18}$. We point out that the closer the global deficiency, $\Delta$, is to zero, the bigger the density of
an arrangement is. It seems like an arrangement in which $\Delta=0$ is not possible.

\subsection{The case ${\cal T}=(3,6,3,6)$}
As before we will use the same technique to find the ``trivial" upper bound which is now usual for these half domination problems.
This upper bound is also sharp in this case.

\begin{figure}[h]
\centering
$\underset{(a)\ Density \ \frac{2}{3} } {\includegraphics[scale=0.2]{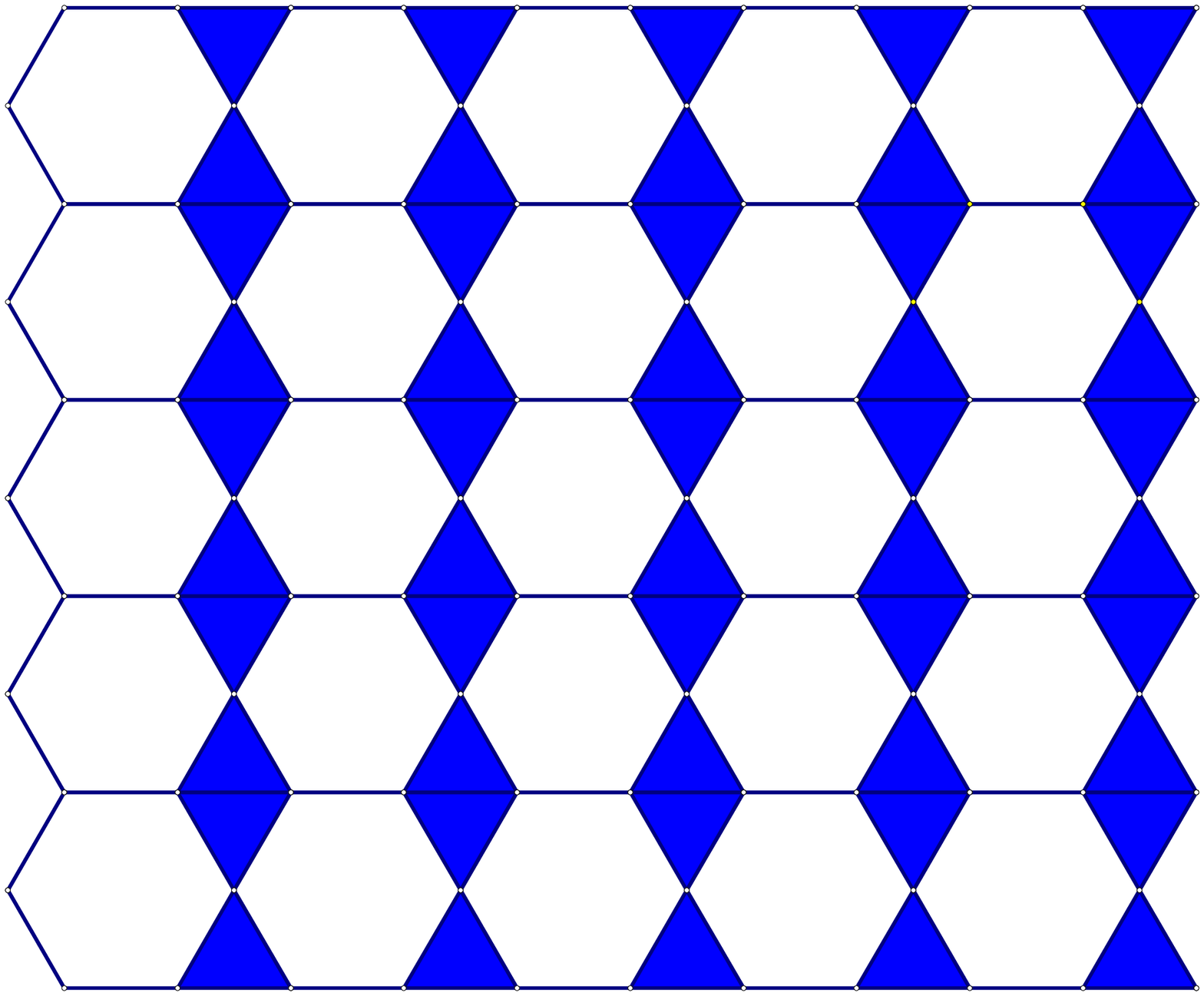}} \ \underset{(b)\ The\ infinity\  column}{\includegraphics[scale=0.4]{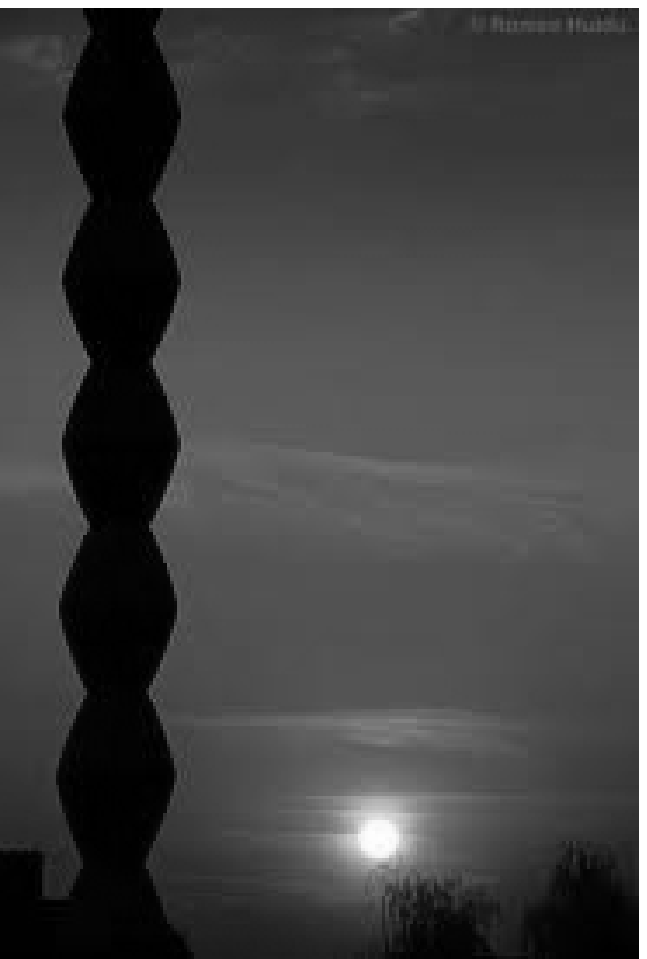}} \ \ \underset{(c)\ Domain}{\includegraphics[scale=0.2]{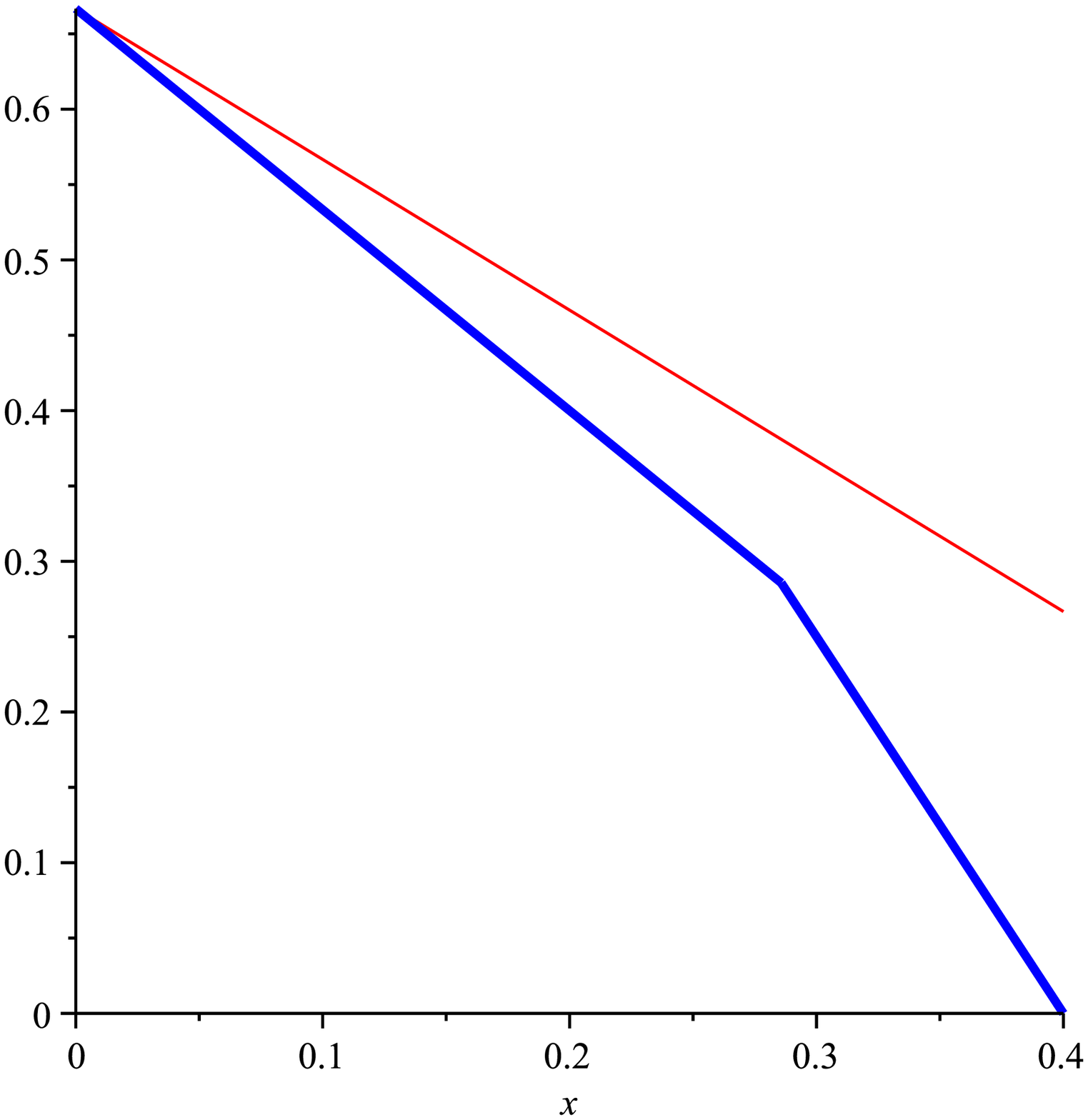} }$
\caption{``Trivial" arrangement} \label{6363arr}
\end{figure}

For a vertex $v$ corresponding to a hexagon we get
$3x_v+x_v^{*}\le 6$, and for $v$ corresponding to a triangle
we have $2x_v+x_v^{*}\le 3$. \n We denote by $H$ the sum of all $x_v$ over all vertices corresponding to hexagons and by $T$ the sum of all $x_v$ over all vertices corresponding to triangles. These inequalities imply

$$3H+(2H+2T)\le 6mn\ \text{and}\ 2T+(4H+T)\le 6mn.$$

\n The system in $x=\frac{H}{3mn}$ and $y=\frac{T}{3mn}$ becomes

$$\begin{cases} 5x+2y\le 2\\ \\
4x+3y\le 2,
\end{cases}$$

\n and the maximum of $x+y$ is equal to $\frac{2}{3}$ and attained for $x=0$ and $y=\frac{2}{3}$ (see Figure~\ref{6363arr} (c)).
Nevertheless, there seems to be a surprising resemblance  between this optimal solution and the well know sculpture of C. {\rm Br\^{\i}ncu\c{s}i}.
(see Figure~\ref{6363arr} (b)).

\subsection{The case ${\cal T}=(3,4,6,4)$}

The system in this case becomes a little more complicated since we have three different type of tilings, but the gives an upper bound which is strictly less than $\frac{2}{3}$.

\begin{theorem}\label{bestupperbound33344}
The half-domination density for the tessellation ${\cal T}=(3,4,6,4)$ satisfies
$$ \rho_{(3,4,6,4)}\le \frac{19}{30}.$$
\end{theorem}

\n \proof. \ The inequalities defining the problem are given by

$$
\begin{cases}
2x_v+x_v^{*}\le 3\ \text{for v corresponding to a triangle}\\ \\
2x_v+x_v^{*}\le 4\ \text{for v corresponding to a square}\\ \\
3x_v+x_v^{*}\le 6\ \text{for v corresponding to a hexagon.}
\end{cases}
$$

As before we introduce $T=\underset{\text{v for triangle}}{ \sum}x_v$, $S=\underset{\text{v for square}}{ \sum}x_v$, and $H=\underset{\text{v for hexagon}}{ \sum}x_v$. The inequalities above give

$$
\begin{cases}
2T+2S\le 3(2mn)\\ \\
2S+(3T+6H)\le 4(3mn)\\ \\
3H+2S\le 6mn,
\end{cases}\ or\
\begin{cases}
2x+2y\le 1\\ \\
3x+2y+6z\le 2\\ \\
2y+3z\le 1,
\end{cases}
$$
\n where $x=\frac{T}{6mn}$, $y=\frac{S}{6mn}$ and $z=\frac{H}{6mn}$. The usual optimization methods give the maximum for $x+y+z$, under the above constrains and $x,y,z\in [0,1)$, to be attained for $x=\frac{1}{5}$, $y=\frac{3}{10}$ and $z=\frac{2}{15}$ and a value of $\frac{19}{30}$.\eproof

\begin{figure}[h]
\centering
\includegraphics[scale=0.5]{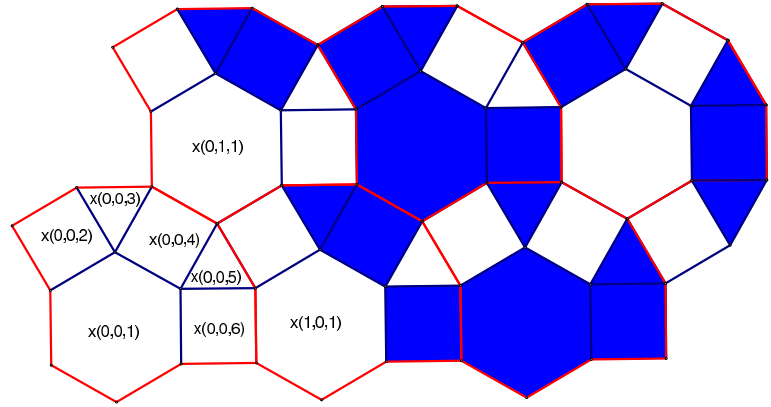}
\caption{Variables convention} \label{6434arrvar}
\end{figure}

The variables we are going to use for writing the system for LPSolve are indexed as before $x_{i,j,k}$, with
$i \in \{0,1,2,...,n-1\}$, $j\in \{0,1,2,...,m-1\}$ and $k\in\{1,2,3,4,5,6\}$ (see Figure~\ref{6434arrvar}).
Given $(i,j)$ as above, we let  $u,v\in \{0,1,2,...,n-1\}$, and $w,t\in \{0,1,2,...,m-1\}$ such that
$u\equiv i-1\ (mod\ n)$, $v\equiv i+1\ (mod \ n)$, $w\equiv j-1\ (mod\ m)$, and  $t\equiv j+1\ (mod \ m)$.
Then the system can be written in the following way (see Figure~\ref{6434arrvar}):

\begin{equation}\label{6434}
\begin{cases}
3x_{i,j,1}+x_{i,j,2}+x_{i,j,4}+x_{i,j,6}+x_{u,j,6}+x_{i,w,4}+x_{v,w,2}\le 6 \\ \\
2x_{i,j,3}+x_{i,j,2}+x_{u,j,4}+x_{u,t,6}\le 3,\ \ 2x_{i,j,5}+x_{i,j,4}+x_{i,j,6}+x_{v,j,2}\le 3,\\ \\
2x_{i,j,2}+x_{i,j,1}+x_{i,j,3}+x_{u,j,5}+x_{u,t,1}\le 4,\ \ 2x_{i,j,4}+x_{i,j,1}+x_{i,j,3}+x_{i,j,5}+x_{v,j,1}\le 4,\ \text{and}\\ \\
2x_{i,j,6}+x_{i,j,1}+x_{i,j,5}+x_{v,j,1}+x_{v,w,3}\le 4.
\end{cases}
\end{equation}

\n Contrary to what is expected LPSolve takes more time to solve these systems even under one hundred variables.
The best densities we have gotten are listed next:

\vspace{0.1in}
\centerline{
\begin{tabular}{|c|c|c|c|c|c|}
  \hline
  (m,n) & (1,1) & (2,2) & (3,3) & (4,4)    \\
  \hline
  $\rho_{(3,4,6,4),(m,n)}$ & $\frac{1}{2}$& $\frac{7}{12}$ & $\frac{31}{54}$ & $\frac{7}{12}$   \\
  \hline
\end{tabular}}
\vspace{0.1in}

\begin{figure}[h]
\centering
$\underset{(a) \ \rho_{(6,4,3,4)}=\frac{3}{4}}{\includegraphics[scale=0.35]{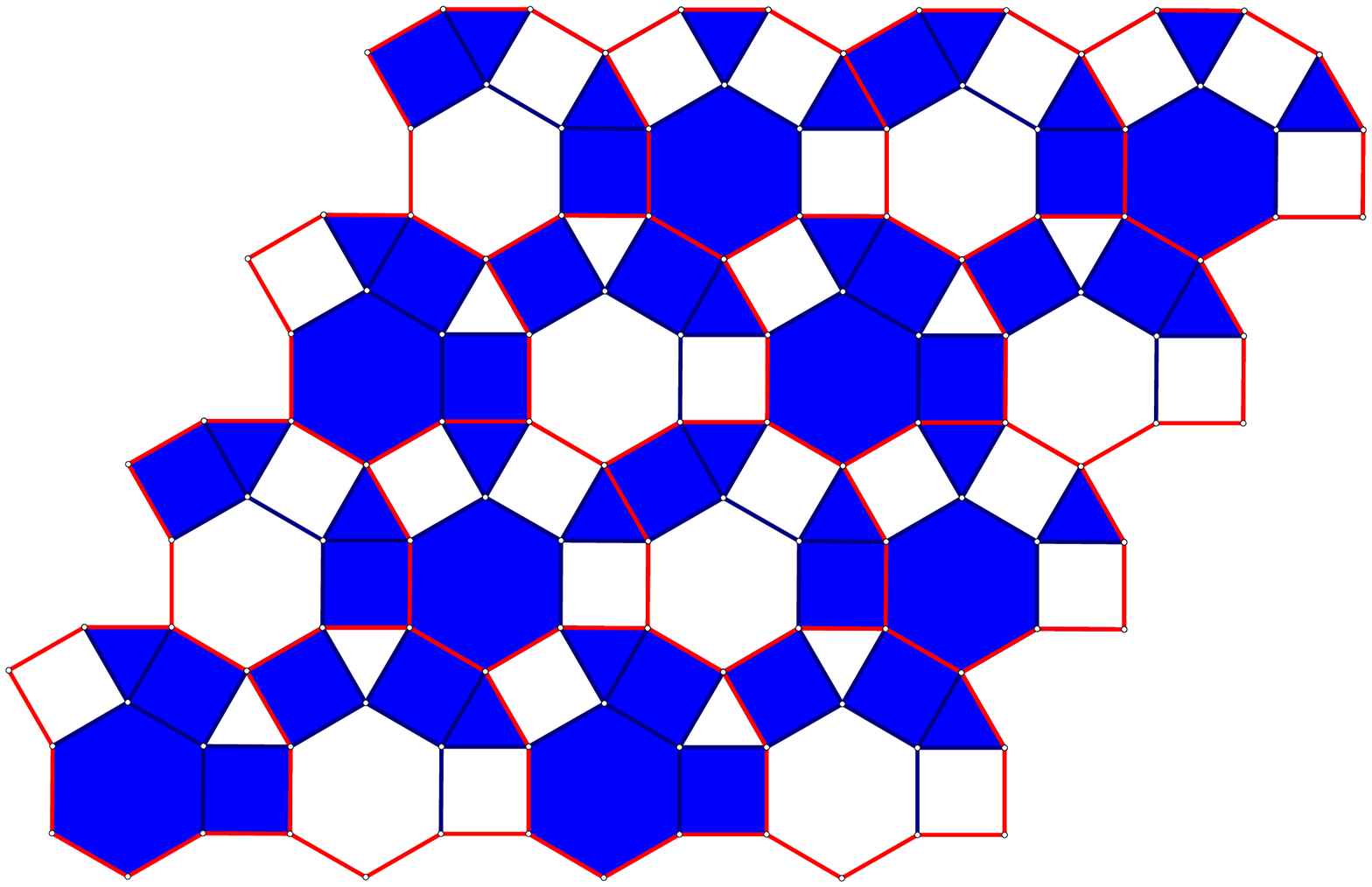}}\ \underset{(b) \ \rho_{(6,3,3,3,3)}=\frac{5}{9}}{\includegraphics[scale=0.4]{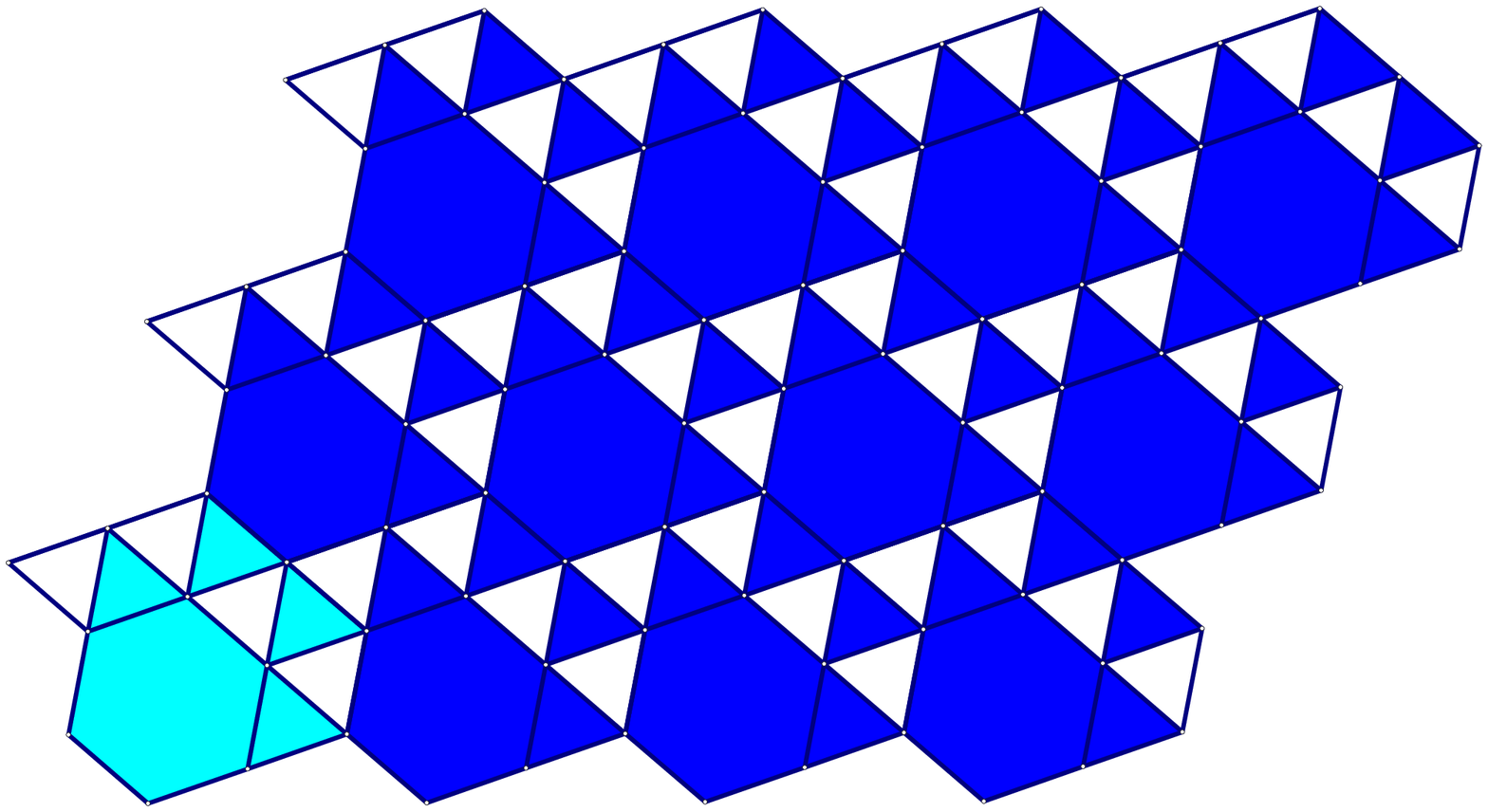}}$
\caption{Some distributions} \label{6434opt}
\end{figure}

\subsection{The rest of the semi-regular tessellations}

\begin{figure}[h]
\centering
$\underset{(a) \ \rho_{(8^2,4)}=\frac{3}{4}}{\includegraphics[scale=0.25]{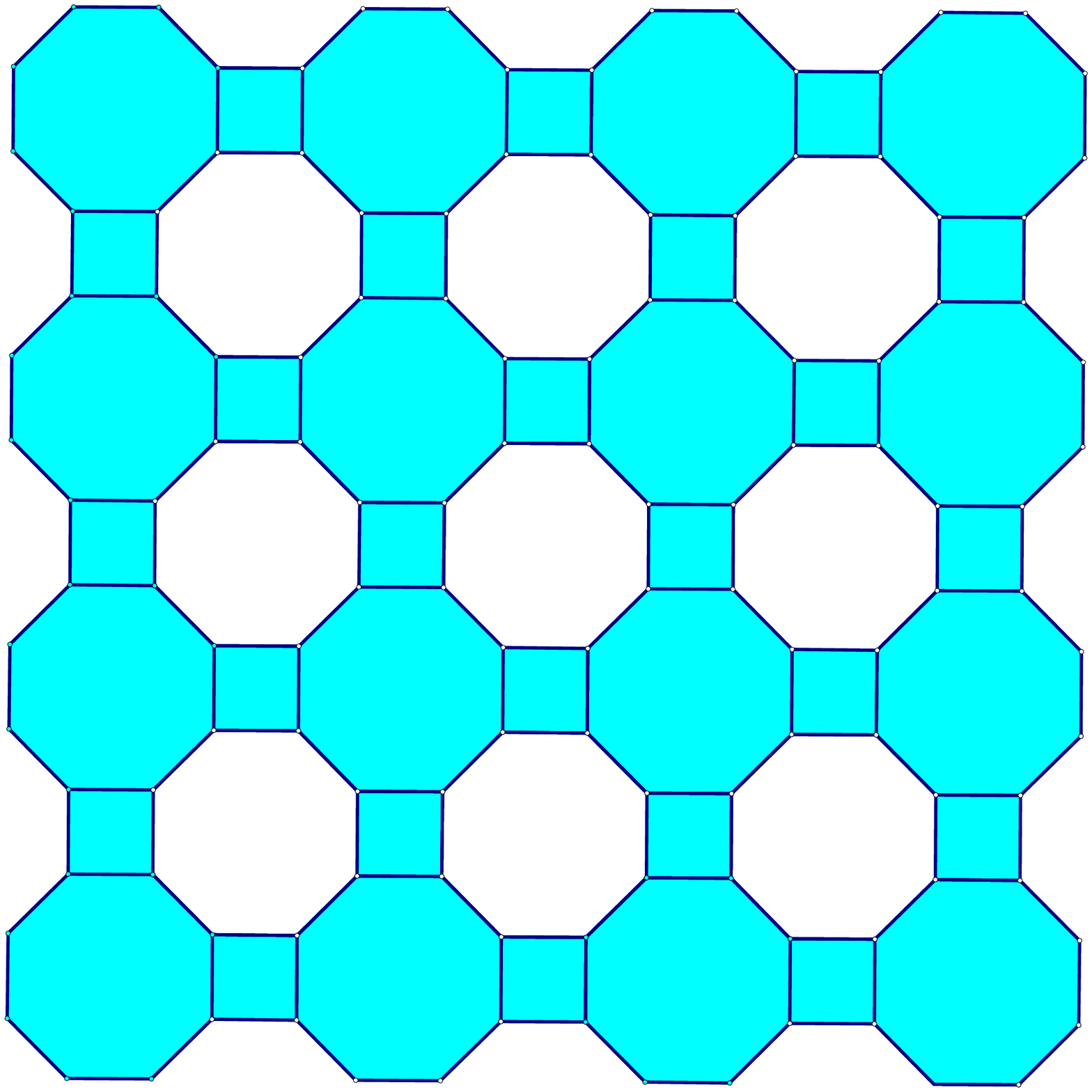}} \ \ \underset{(b) \ \rho_{(12,6,4)}=\frac{5}{6}}{\includegraphics[scale=0.3]{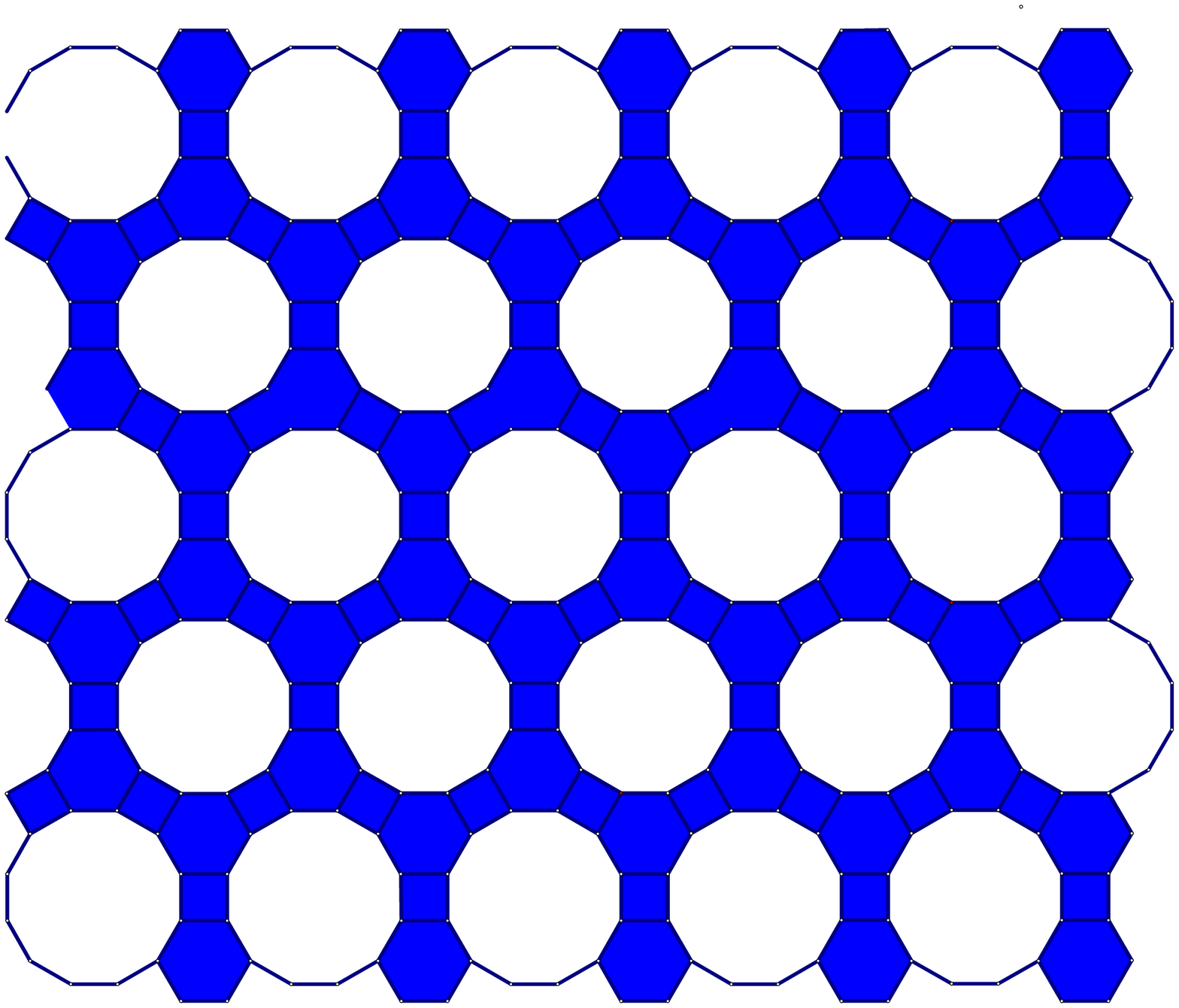}}$
\caption{Best arrangements } \label{884opt}
\end{figure}

We are going to include for the rest of the semi-regular tessellations  only the most significant
facts found but without proves. One case use the same methods to check them.
For the semi-regular tessellations $(8^2,4)$ and $(12,6,4)$ the following arrangements (Figure~\ref{884opt}(a) and (b)) gives the best
densities.

We notice that the deficiency for each of the arrangements in Figure~\ref{884opt}(a) and (b) is equal to zero.

\begin{figure}[h]
\centering
$\underset{(a) \ \rho_{(12^2,3)}=\frac{3}{4}}{\includegraphics[scale=0.25]{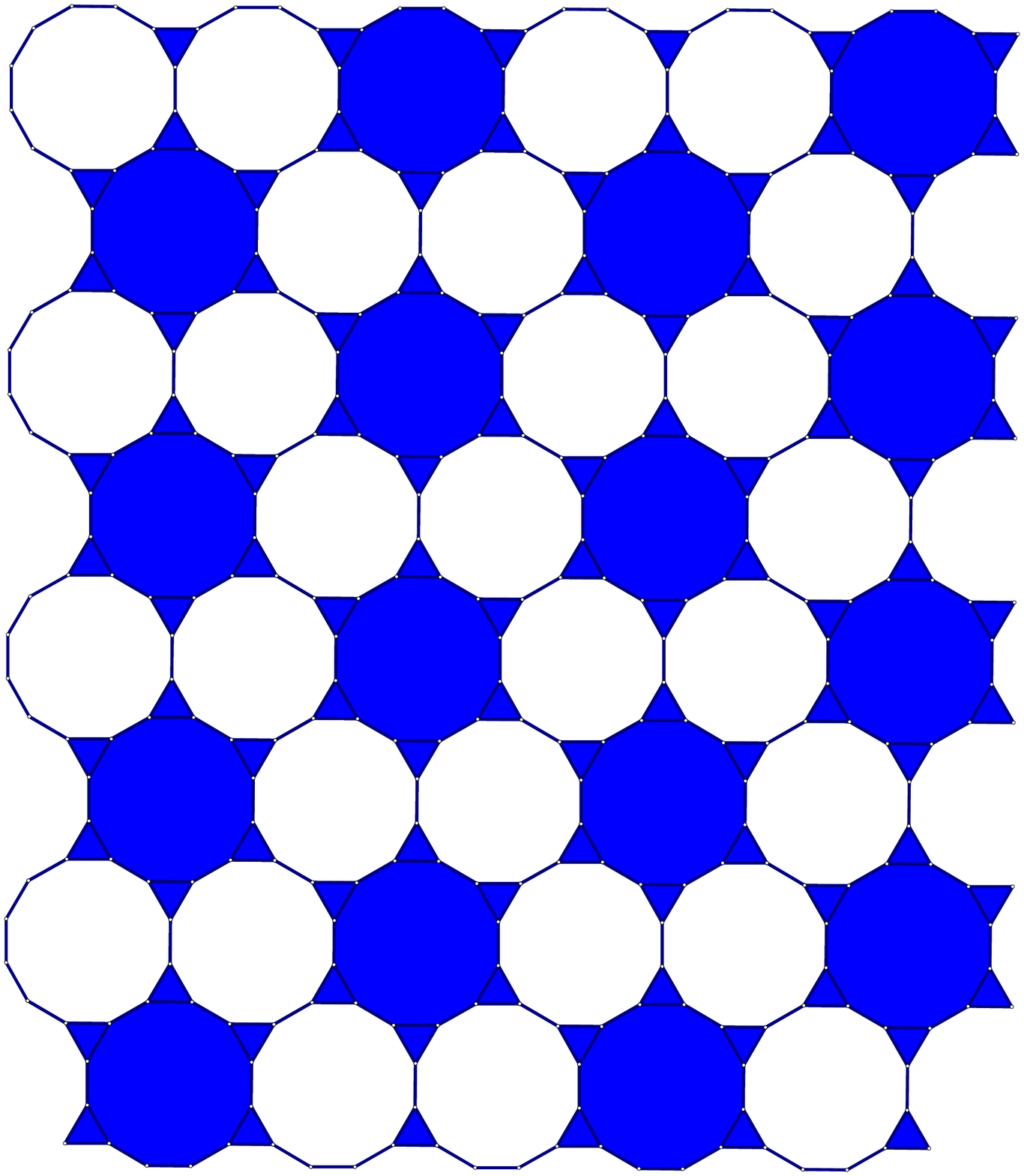}} \ \ \underset{(b) \ \rho_{(4,3,3,4,3)}=\frac{2}{3}}{\includegraphics[scale=0.3]{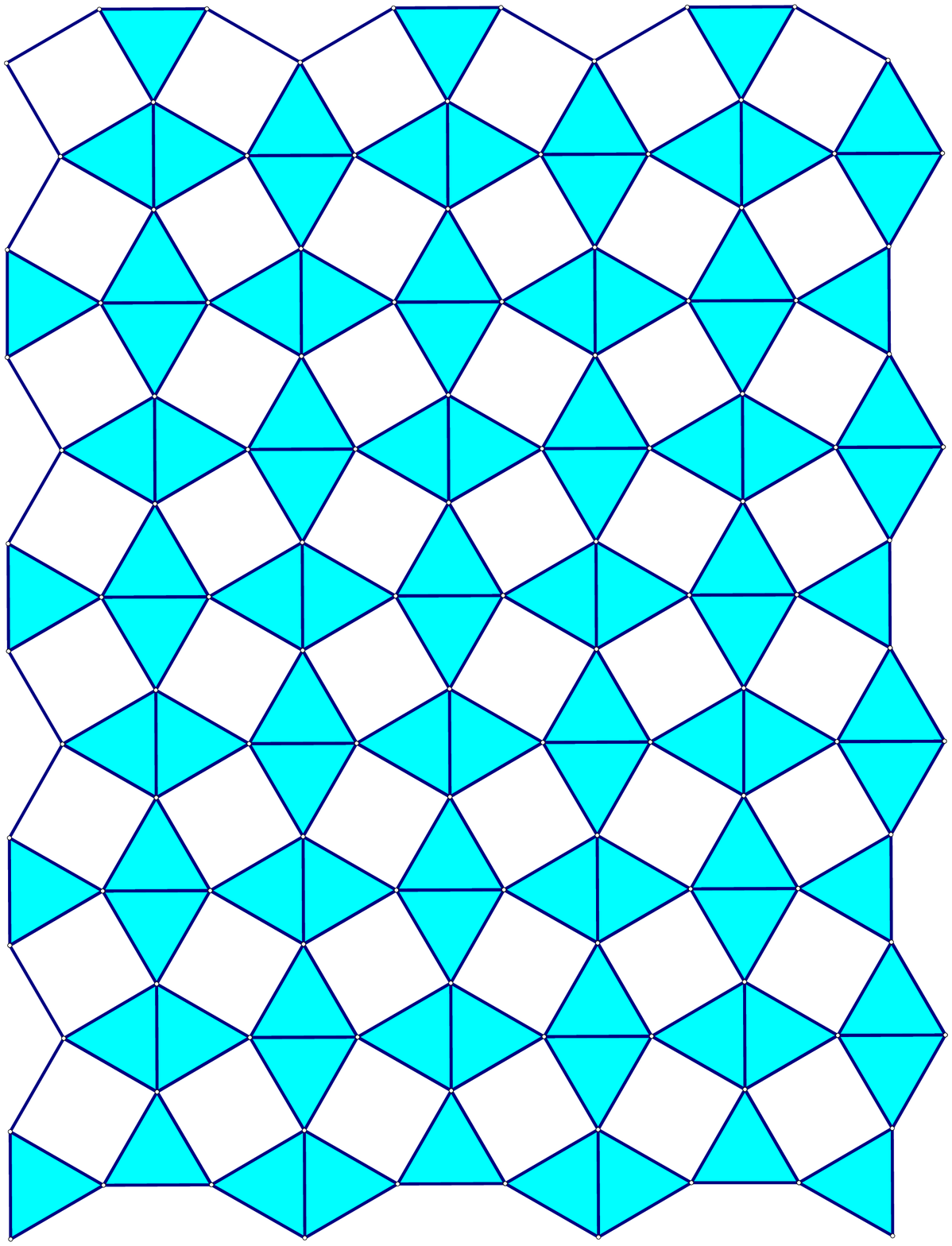}}$
\caption{Best arrangements } \label{12312opt}
\end{figure}

For the semi-regular tessellation $(12^2,3)$, we have gotten the theoretical upper bound of $\frac{7}{9}$ and one of the arrangements as in Figure~\ref{12312opt}(a)
which we think it is actually the sharp. In Figure~\ref{12312opt}(b), we have a sharp arrangement for tessellation $(4,3,3,4,3)$.
Finally, for $\cal T=(6,3,3,3,3)$ we see see an arrangement of density $\frac{5}{9}$ in Figure~\ref{6434opt} (b), but one can show that the best density is
actually $\frac{2}{3}$ given by the distribution shown in Figure~\ref{3^46mn84}.

\section{Conjectures and other comments}

From what we have seen so far, there are some patterns that emerge. Given a {\it vertex transitive}  infinite graph (for every two vertices, there exists a graph isomorphism mapping one
vertex into the other) have half-domination arrangements which have rational best densities. It is not clear if such arrangements are {\it unique} (up to the isomorphisms of the graph)  or there exist essentially different variations. In any case, we see that if the deficiency is zero, then the solution seems to be unique. If the deficiency is positive, one may expect to have more solutions and we have such an example in the case of the King's Graph (see \cite{ipw}). The bigger the deficiency the higher the number of combinatorial possibilities that can result in maximum arrangements but, we conjecture that there are only finitely many of them. Results that show the exact number of such maximum arrangements are, nevertheless, at our interest in further investigations. However, we believe that the right methods to approach these questions successfully, even with the assistance of powerful computers, are yet to be discovered.
Another path of investigations is to look into finding similar answers to $k$-dependence problems in all of the graphs studied here.

\noindent {\em Math Department, Columbus State University, 4225 University Avenue,
Columbus, GA 31907, e-mail: ionascu@columbusstate.edu}

\end{document}